\documentclass[10pt]{elsarticle}
\usepackage{amsmath}
\usepackage{amssymb}
\usepackage{graphicx,epsfig}
\usepackage{latexsym}
\newtheorem{prop}{Proposition}
\newtheorem{lemma}{Lemma}

\newtheorem{theorem}{Theorem}
\newtheorem{remark}{Remark}

\newcommand{\R}{\mathbb{R}}
\newcommand{\N}{\mathbb{N}}

\font\tenmath=msbm10 \font\sevenmath=msbm7 \font\fivemath=msbm5
\newfam\mathfam \textfont\mathfam=\tenmath
\scriptfont\mathfam=\sevenmath \scriptscriptfont\mathfam=\fivemath
\def\math{\fam\mathfam}

\def \={{\buildrel {\rm (law)} \over =}}
\def \N{{\math N}}
\def \R{{\math R}}

\def\L{\Lambda}

\def\L{\mathbb{L}}
\def\R{\mathbb{R}}

\def\P{\mathrm P}

\def\text{\mbox}

\def\1{{\bf 1}}
\def\qed{ \hfill \vrule width.25cm height.25cm depth0cm\smallskip}

\newcommand{\basa}{\begin{assumption}}
\newcommand{\easa}{\end{assumption}}

\newcommand{\Cov}{\mbox{Cov}\,}
\newcommand{\Var}{\mbox{Var}\,}

\def\limiteN{\renewcommand{\arraystretch}{0.5}
\begin{array}[t]{c}\stackrel{}{\longrightarrow} \\
{\scriptstyle N\rightarrow
\infty}\end{array}\renewcommand{\arraystretch}{1}}

\def\limiteNa{\renewcommand{\arraystretch}{0.5}
\begin{array}[t]{c}\stackrel{}{\longrightarrow} \\
{\scriptstyle N_a\rightarrow
\infty}\end{array}\renewcommand{\arraystretch}{1}}

\def\limiteloiN{\renewcommand{\arraystretch}{0.5}
\begin{array}[t]{c}\stackrel{{\cal D}}{\longrightarrow} \\
{\scriptstyle N\rightarrow
\infty}\end{array}\renewcommand{\arraystretch}{1}}

\def\limiteloin{\renewcommand{\arraystretch}{0.5}
\begin{array}[t]{c}\stackrel{{\cal D}}{\longrightarrow} \\
{\scriptstyle n\rightarrow
\infty}\end{array}\renewcommand{\arraystretch}{1}}

\def\limiteloiNa{\renewcommand{\arraystretch}{0.5}
\begin{array}[t]{c}\stackrel{{\cal D}}{\longrightarrow} \\
{\scriptstyle N_{a}\rightarrow
\infty}\end{array}\renewcommand{\arraystretch}{1}}

\def\limiteL2{\renewcommand{\arraystretch}{0.5}
\begin{array}[t]{c}\stackrel{\L^{2}(\Omega)}{\longrightarrow} \\
{\scriptstyle N\rightarrow
\infty}\end{array}\renewcommand{\arraystretch}{1}}

\def\simloi{\renewcommand{\arraystretch}{0.5}
\begin{array}[t]{c}\stackrel{{\cal D}}{\sim} \\
{}\end{array}\renewcommand{\arraystretch}{1}}

\def\limiteprobaN{\renewcommand{\arraystretch}{0.5}
\begin{array}[t]{c}\stackrel{\P}{\longrightarrow} \\
{\scriptstyle N_{a}\rightarrow
\infty}\end{array}\renewcommand{\arraystretch}{1}}

\begin{document}

\title{A wavelet analysis of the Rosenblatt process: chaos expansion and estimation  of the self-similarity parameter}
\author[sam]{J.-M. Bardet\corref{cor1}}
\ead{bardet@univ-paris1.fr}
\author[lpp]{C.A. Tudor $\dag ,$ \footnote{ $\dag$ Associate member: SAMM,  Universit\'{e} de Paris 1, 90, rue de Tolbiac, 75634, Paris, France}}
\ead{tudor@math.univ-lille1.fr}
\cortext[cor1]{Corresponding author}
\address[sam]{SAMM,  Universit\'{e} de Paris 1, 90, rue de Tolbiac, 75634, Paris, France}
\address[lpp]{Laboratoire Paul Painlev\'e, Universit\'e de Lille 1,
 F-59655 Villeneuve d'Ascq, France.}
\begin{abstract}
By using chaos expansion into multiple stochastic integrals, we make a wavelet analysis of two self-similar stochastic processes: the fractional Brownian motion and the Rosenblatt process. We study the asymptotic behavior of the  statistic based on the wavelet coefficients of these processes.
Basically, when applied to a non-Gaussian process (such as the Rosenblatt process) this statistic satisfies a non-central limit theorem even when we increase the number of vanishing moments of the wavelet function. We apply our limit theorems to construct  estimators for the self-similarity index and we illustrate our results by  simulations.
\end{abstract}

\maketitle

\section{Introduction}
The self-similarity property for a stochastic process means that scaling of time is equivalent to an appropriate scaling of
space. That is, a process $(Y_{t})_{t\geq 0}$ is self-similar of order $H>0$ if for all $c>0$ the processes $(Y_{ct}) _{t \geq
0} $ and $(c^{H} Y _{t} )_{t\geq 0}$  have the same finite dimensional distributions. This property is crucial in applications
such as network traffic analysis, mathematical finance, astrophysics, hydrology  or image processing. We refer to the
monographs \cite{Beran}, \cite{EM} or \cite{ST} for complete expositions on theoretical and practical aspects of self-similar stochastic  processes.

The most popular self-similar process is the fractional Brownian motion (fBm).  Its practical applications are numerous.  This
process  is defined as a centered Gaussian process $(B^{H}_{t}) _{t\geq 0}$ with covariance function
\begin{equation}\label{covFBM}
R_H(t,s)=\mathbb{E} (B^{H}_{t}B^{H} _{s}) =\frac{1}{2} \left( t^{2H } + s^{2H} -\vert t-s \vert ^{2H}\right), \hskip0.5cm t,s \geq 0.
\end{equation}
It can be also  defined as the only Gaussian self-similar process with stationary increments. Recently, this stochastic process
has been widely studied from the stochastic calculus point of view as well as from the  statistical analysis point of view.
Various types of stochastic integrals with respect to it have been introduced  and several  types of stochastic differential
equations driven by fBm have been considered. Other  stochastic processes which are self-similar with stationary increments
are the Hermite processes (see \cite{BrMa}, \cite{DM}, \cite{Ta1}); an Hermite process of order $q$  is  actually an iterated
integral of a deterministic function with  $q$ variables with respect to the standard Brownian motion.  These processes
appear as limits in the so-called noncentral limit theorem and they have the same covariance as the fBm. The fBm is obtained
for $q=1$ and it is the only Gaussian Hermite process. For $q=2$ the  corresponding  process is known as the Rosenblatt
process. Although it received a less important attention than the fractional Brownian motion, this process is still of
interest in practical applications because of its self-similarity, stationarity and long-range dependence of increments.
Actually the great popularity  of the fractional Brownian motion in practice (hydrology, telecommunications) is due to
these properties; one prefers fBm rather than higher order Hermite
process because it is a Gaussian process  and the corresponding stochastic analysis is much easier. But
in concrete situations when empirical data attest to the presence of
self-similarity and long memory without the Gaussian property (an example is mentioned in the paper \cite{Taqqu3}), one can use a
Hermite process living in a higher chaos, in particular the Rosenblatt process.

When studying self-similar processes,  a question of major interest is to estimate their self-similarity order. This is
important because the self-similarity order characterizes in some sense the process: for example in the fBm case as well as for
Hermite processes this order gives the long-range dependence property of its increments and it characterizes the regularity of the trajectories. Several statistics, applied directly to the process or to its increments, have been introduced to address  the problem of estimating the self-similarity index. Naturally, parametric statistics (exact or  Whittle approached maximum likelihood) estimators have been considered. But to enlarge the method to a more general class of models (such as e.g. locally or asymptotically self-similar models), it could be interesting to apply semiparametric methods such as  wavelets based, log-variogram or log-periodogram estimators. Information and details on these
various approaches can be found in the books of Beran \cite{Beran} and Doukhan {\it et al.} \cite{Doukhan}.

Our purpose is to develop a wavelet-based analysis of the fBm and the  Rosenblatt process by using multiple Wiener-It\^o integrals and to apply our results to estimate the self-similarity index. More precisely, let $\psi: \mathbb{R}\to \mathbb{R}$ be a
continuous function with support included in the interval $[0,1]$ (called "mother wavelet"). Assume that there exists an
integer $Q\geq 1$ such that
\begin{equation}
\label{cond-psi} \int_{\mathbb{R}} t^{p} \psi (t) dt =0 \mbox{ for } p=0, 1, \ldots , Q-1
\end{equation}
and \begin{equation*} \int_{\mathbb{R} }   t^{Q} \psi(t)dt \neq 0.
\end{equation*}
The integer $Q\geq 1$ is called the {\it number of vanishing moments}. For a stochastic process $(X_{t}) _{t\in [0,N] } $ and for a "scale" $a\in \mathbb{N}^{\ast}$ we define its wavelet
coefficient by
\begin{equation} \label{coef}
d(a,i)= \frac{1}{\sqrt{a}} \int_{-\infty } ^{\infty } \psi \left( \frac{t}{a}-i\right) X_{t}dt = \sqrt{a} \int_{0}^{1} \psi(x)
X_{a(x+i) } dx
\end{equation}
for $i=1,2, \ldots ,N_a$ with $N_a= [N/a] -1$. Let us also introduce the normalized wavelet coefficient
\begin{equation}\label{coefnorm}
\tilde{d}(a,i)  =\frac{ d(a,i)}{\left( \mathbb{E} d^2(a,i)\right) ^{\frac{1}{2}}}
\end{equation}
and the statistic
\begin{equation} \label{VN}
 V_{N}(a)= \frac{1}{N_{a} }\sum_{i=1}^{N_{a}} \left( \tilde{d} ^{2} (a,i) -1\right) .
\end{equation}
The wavelet analysis consists in the study of  the asymptotic  behavior of the sequence $V_{N}(a)$ when $N\to \infty$. Usually, if $X$ is  a stationary long-memory process or a self-similar second-order process, then $\mathbb{E} d^2(a,i)$ is a power-law function of $a$ with  exponent $2H-1$ (when $a \to \infty$) or $2H+1$ respectively. Therefore, if $V_N(a)$  converges to $0$, a log-log-regression of $\frac{1}{N_{a} }\sum_{i=1}^{N_{a}} d^2  (a_j,i)$ onto $a_j$ will provide an estimator of $H$ (with an appropriate choice of $(a_j)_j$). Hence, the asymptotic behavior of $V_{N}(a)$ will completely give the asymptotic behavior of the estimator (see the Section \ref{Appli} for
details). There are four main advantages to use such an estimator based on wavelets: firstly, it is  a semiparametric method that can be easily generalized. Secondly, it is based on the log-regression of the wavelet coefficient sample variances onto several scales and the graph of such a regression provides interesting information concerning the goodness-of-fit of the model ($\chi^2$ goodness-of-fit have been  defined and studied in \cite{Bardet} or \cite{BB} from this log-regression). Thirdly, it is a computationally efficient  estimator (this is due to the Mallat's algorithm for computing the wavelet coefficients). Finally, it is a very robust method: it is not sensitive to possible
polynomial trends as soon as the number of vanishing moments $Q$ is large enough.

 This  method has been introduced by Flandrin \cite{flan2} in the case of fBm and later it has been extended to more general processes in  \cite{Abry2}. The asymptotic behavior of the estimator obtained from the wavelet variation has been  specified in the case of long-memory Gaussian processes in \cite{BLMS} or \cite{MRT2}, for long memory linear processes in \cite{RT} and for  multiscale fractional Gaussian processes in \cite{BB}. However the case of the Rosenblatt process has not yet been studied (note that the wavelet synthesis of Rosenblatt processes was treated in \cite{pipiras} and we will use this method in the section devoted to simulations).
 
 In our paper we will use a recently developed theory
based on Malliavin calculus and  Wiener-It\^o multiple stochastic integrals. Let us briefly recall these new results. In
\cite{NP} the authors gave necessary and sufficient conditions for a sequence of random variables in a fixed Wiener chaos
(that means in essence that these random variables are iterated integrals of a fixed order with respect to a given Brownian
motion) to converge to a standard Gaussian random variable (one of these conditions is that the sequence of the fourth order moments
converges to $3$ which represents the moment of order $4$ of a standard Gaussian random variable). Another equivalent condition is
given in the paper \cite{NOT} in terms of the Malliavin derivative. These results created a powerful link between the
Malliavin calculus and limit theorems and they have already been used in several papers (for example in \cite{TV},
\cite{CNT} and \cite{CNT1} to study the variations of the Hermite processes).\\

Recall (see \cite{flan2} and \cite{Bardet}) that if $X=B^{H}$ is a fBm in (\ref{VN}) then the following fact happens: for any $Q>1$ and $H\in (0,1)$ the
statistic $V_{N}(a)$ renormalized by $\sqrt{N}$ converges to a centered Gaussian random variable. If $Q=1$ then the barrier
$H=3/4$ appears: that is, the behavior of $V_{N}(a)$ is normal ({\it i.e.} it satisfies a central limit theorem) only if
$H\in (0,3/4) $ and we prove in Section \ref{SecFBM} that the limit of $V_{N}(a)$ (normalized by $N^{2-2H}$) is a
Rosenblatt random variable $R^{H}_{1}$ when $H\in (3/4, 1)$. This non-central limit theorem is new in the case of the wavelet based statistic (it is known for other type of statistics, such as the statistic constructed from the variations of the process, see \cite{TV}). In this case we also prove that the limit holds only  in law and not in $\L^2$ in contrast to  the case of quadratic variations studied in \cite{TV}.\\

The study of $V_{N}(a)$ in the Rosenblatt process case (see Section \ref{SecFBM}, formula (\ref{repZ}) for the definition of the Rosenblatt process) with $H\in (1/2,1)$ put in light interesting and somehow intriguing phenomena. {\em The main fact is that
the number of vanishing moments $Q$ does not affect its convergence and the limit of $V_{N}(a)$ is always non-Gaussian }(it is still
Rosenblatt). This fact is unexpected and different to the situations known in the literature. Actually, the statistic $V_N$ can be decomposed into two parts: a term in the fourth chaos (an iterated
integral of order $4$ with respect to a Wiener process) and a term in the second chaos. We analyze here both terms and we
deduce that the term in the fourth Wiener chaos keeps some of the characteristics of the Gaussian case (it has to be
renormalized by $\sqrt{N}$ and it has a Gaussian limit for $H\in (1/2,3/4) $). But the main term in $V_{N}(a)$ which gives
the normalization  is the second chaos term and its detailed analysis shows that the normalization depends on $H$ (it is of order of $N^{1-H}$) and its limit
is (in law) a Rosenblatt random variable. \\

The consequences of these results are also interesting for the wavelet based estimator of the self-similarity order. Assume that a sample  $(X_1,\cdots,X_N)$ is obseved, where $X$ is a fBm or a Rosenblatt process. First, by approximating the wavelet coefficients (\ref{coef}), we consider a statistic $\widehat V_N(a)$ that can be computed from the observations  $(X_1,\cdots,X_N)$ and we prove that the limit theorems satisfied by $V_N(a)$ also hold for  $\widehat V_N(a)$ as soon as $a$ is large enough with respect to $N$. Secondly, we deduce the convergence rates for the wavelet based estimator of $H$ following the cases: $Q\geq 2$ and $X$ is a fBm, $Q=1$, $H\in (3/4,1)$ and $X$ is a fBm or $H\in (\frac{1}{2}, 1)$ and  $X$ is a Rosenblatt process. The regularity of $\psi$ also plays  a role. For practical use, it is clear that if $X$ is a fBm it is better to chose $Q\geq 2$ and the mother wavelet $\psi$ to be twice continuously differentiable (for example, this is the case when $\psi$ is a Daubechies wavelet with order $\geq 8$). On the other hand, if  $X$ is a Rosenblatt process the number of vanishing moments  $Q$ plays no role. Our simulations illustrate the convergence of $\widehat V_N(a)$ and of the estimator of $H$ in this last case.

 Our results open other related questions: for a stationary long-memory Gaussian or linear process,  parametric estimators (such as the Whittle's estimator, see  \cite{FT} and \cite{GS}) or  semiparametric estimators (such as wavelet based estimator, see  \cite{MRT2} and \cite{RT}) satisfy a central limit theorem with an  ``usual'' convergence rate (that is $\sqrt N$ for Whittle's estimator and $\sqrt{N/a}$ for wavelet based estimator). This is not the case of the Whittle's estimator for non-linear functionals of a long-memory Gaussian process (see \cite{GT}). As we prove in our paper,  this is no longer the case for the wavelet  estimator based on the observation of a  Rosenblatt process. 
Therefore it seems that the  Rosenblatt processes and the second order polynomials of Gaussian process give the  same asymptotic behavior for the  estimators. Indeed, in \cite{GT} it has been  established that  the Whittle's estimator in the case of second order polynomials of Gaussian processes satisfies a noncentral limit theorem with a Rosenblatt distribution as a limit and with the same convergence rate as in Theorem \ref{Multros} below.  Several interesting questions arise naturally from this: firstly, are there any  reciprocal results true? ({\it i.e.}  may we suspect that the Whittle's estimator of the long-memory parameter for the increments of the Rosenblatt process and the wavelet based estimator of the long-memory parameter for a second order polynomial of a long-memory Gaussian process also satisfy the same noncentral limit theorem?). Secondly, are the conclusions the same for other semiparametric estimators such as log-periodogram or local Whittle's estimators?    \\

We organized the paper as follows. Section \ref{Prem} contains some preliminaries on multiple Wiener-It\^o integrals with respect to the Brownian motion. In Section \ref{SecFBM} we treat the situation when the driven process is the fBm. In this case our new result is the noncentral limit theorem  satisfied by the wavelet based statistic proved in Theorem 2. In Section \ref{SecROS} we enter into a non-Gaussian world: our observed process is the Rosenblatt process and using the techniques of the Malliavin calculus and recent interesting results for the convergence of sequence of multiple stochastic integrals, we study in details the sequence $V_{N}(a)$. In Section \ref{Appli}  we construct an observable estimator based on the approximated wavelet coefficients and we study its asymptotic behavior. We also construct an estimator of the self-similarity order and we illustrate its convergence by numerical results. 
\section{Preliminaries}\label{Prem}

\subsection{Basic tools on multiple Wiener-It\^o integrals}

Let $(W_{t})_{t\in [0,T] }$ be a classical Wiener process on a standard Wiener space $\left( \Omega , {\mathcal{F}}, P\right)
$. If $f\in \L^{2}([0,T]^{n}) $ with $n\geq 1$ integer, we introduce the multiple Wiener-It\^o integral of $f$ with respect to
$W$. We refer to \cite{N} for a detailed exposition of the construction and the properties of multiple Wiener-It\^o integrals.\\
~\\
Let $f\in {\mathcal{S}}_{n}$, that means that there exists $n \geq 1$ integer such that
\begin{equation*}
f=\sum_{i_{1}, \ldots , i_{n}}c_{i_{1}, \ldots ,i_{n} }1_{A_{i_{i}}\times \ldots \times A_{i_{n}}}
\end{equation*}
where the coefficients satisfy $c_{i_{1}, \ldots ,i_{n} }=0$ if two indices $i_{k}$ and $i_{\ell}$ are equal and the sets
 $A_{i}\in {\mathcal{B}}([0,T]) $ are disjoint. For such a step function $f$ we define

\begin{equation*}
I_{n}(f)=\sum_{i_{1}, \ldots , i_{n}}c_{i_{1}, \ldots ,i_{n} }W(A_{i_{1}})\ldots W(A_{i_{n}})
\end{equation*}
where we set $W([a,b])= W_{b}-W_{a}$. It can be seen that the application $I_{n}$ constructed above from ${
\mathcal{S}}_{n}$ equipped with the scaled norm $\frac {1}{\sqrt{n!}}\Vert\cdot\Vert_{\L^{2}([0, T] ^{n})}$ to $\L^{2}(\Omega)$ is
an isometry on ${\mathcal{S}}_{n}$, i.e. for $m,n$ positive integers,
\begin{eqnarray}
\mathbb{E}\left(I_{n}(f) I_{m}(g) \right) &=& n! \langle f,g\rangle _{\L^{2} ([0,T]^{n})}\quad \mbox{if } m=n, \label{isom}\\
\mathbb{E}\left(I_{n}(f) I_{m}(g) \right) &= & 0\quad \mbox{if } m\not=n. \nonumber
\end{eqnarray}
It also holds that
\begin{equation}\label{ftilde}
I_{n}(f) = I_{n}\big( \tilde{f}\big)
\end{equation}
where $\tilde{f} $ denotes the symmetrization of $f$ defined by $\tilde{f}%
(x_{1}, \ldots , x_{n}) =\frac{1}{n!} \sum_{\sigma \in {\cal S}_{n}} f(x_{\sigma (1) }, \ldots , x_{\sigma (n) } ) $.\\
Since the set ${\mathcal{S}}_{n}$ is dense in $\L^{2}([0,T]^{n})$ for every $n\geq  1,$ the mapping $%
I_{n}$ can be extended to an isometry from $\L^{2}([0,T]^{n}) $ to  $%
\L^{2}(\Omega)$ and the above properties (\ref{isom}) and (\ref{ftilde}) hold true for this extension. Note also that $I_{n}$ can be viewed as an iterated
stochastic integral
\begin{equation*}
I_{n}(f)=n!\int_{0}^{1}\int_{0} ^{t_{n}} \ldots \int_{0}^{t_{2}}f(t_{1}, \ldots , t_{n}) dW_{t_{1}}\ldots dW_{t_{n}}.
\end{equation*}
We recall the product formula  for two multiple integrals (see \cite{N}): if $f\in \L^{2}([0,T] ^{n})$ and $g\in \L^{2}([0,T]  ^{m} )$ are symmetric functions, then it holds that
\begin{equation}  \label{prod}
I_{n}(f)I_{m}(g)= \sum_{\ell=0}^{m\wedge n} \ell! C_{m}^{\ell}C_{n}^{\ell} I_{m+n-2\ell} (f\otimes _{\ell}g)
\end{equation}
where $C_{m}^{\ell}=\frac{m!}{\ell ! (m-\ell)!} $ and  the contraction $f\otimes _{\ell}g$ belongs to $\L^{2}([0,T]^{m+n-2\ell}) $ for $\ell=0, 1, \ldots , m\wedge n$ with
\begin{eqnarray}\nonumber
(f\otimes _{\ell} g) ( s_{1}, \ldots , s_{n-\ell}, t_{1}, \ldots , t_{m-\ell})  =\int_{[0,T] ^{\ell} } f( s_{1}, \ldots , s_{n-\ell}, u_{1}, \ldots ,u_{\ell})g(t_{1}, \ldots , t_{m-\ell},u_{1},
\ldots ,u_{\ell}) du_{1}\ldots du_{\ell}. \label{contra}
\end{eqnarray}
In general the contraction $f\otimes _{\ell} g $ is not necessarily  a symmetric function even if the two functions $f$ and $g$ are symmetric.

\section{The case of the fractional Brownian motion}\label{SecFBM}
\subsection{A presentation using chaos expansion}
We will assume in this part that $X=B^{H}$ a (normalized) fractional Brownian motion (fBm in the sequel) with Hurst parameter $H\in (0,1)$. Recall that  $B^{H}$ is
a centered Gaussian process with covariance function (\ref{covFBM}). It is the only normalized Gaussian $H$-self-similar process with stationary increments. Recall also the fBm
$(B^{H}_{t})_{t\in [0,N]}$ with Hurst parameter $H\in (0,1)$ can be written
$$B_t^H=\int_0^t K^{H}(t,s)dW_s, \quad t\in[0,N]$$
where $(W_t,t\in[0,N])$ is a standard Wiener process and for $s< t$, and  for $H>\frac{1}{2}$ the kernel $K^{H}$  has the expression
\begin{equation}\label{defK}K^H(t,s)=c_Hs^{\frac{1}{2}-H}\int_s^t(u-s)^{H-\frac{3}{2}}u^{H-\frac{1}{2}}du\end{equation}
with $c_H=\left(\frac{H(2H-1)}{\beta(2-2H,H-\frac{1}{2})}\right)^{\frac{1}{2}}$ and $\beta(\cdot,\cdot)$ the
beta function. For $t>s$ and for every $H\in (0,1)$ the partial derivative of $K^{H}$ with respect to its first variable is given by
\begin{equation}
\label{dK} \frac{\partial K^H}{\partial t} (t,s)= \partial_1 K^H (t,s)=c_{H}\left( \frac{s}{t} \right) ^{\frac{1}{2}-H} (t-s)^{H-\frac{3}{2}}.
\end{equation}
In this case it is trivial to decompose in chaos the wavelet coefficient $d(a,i)$. By a stochastic Fubini theorem we can write
 \begin{eqnarray*}
d(a,i)&=& \sqrt{a} \int_{0}^{1} \psi (x) B^{H} _{a(x+i)} dx =\sqrt{a} \int_{0}^{1} \psi(x) dx \left( \int_{0} ^{a(x+i) }
dB^{H}_{u}\right) \\
&=& \sqrt{a}  \int_{0}^{1} \psi (x) dx \int_{0}^{a(x+i)}K^{H} (a(x+i), u)dW_{u} = I_{1} \left(   f_{a,i} (\cdot )\right)
\end{eqnarray*}
where $I_{1} $ denotes the multiple integral of order one (actually, it is the Wiener integral with respect to $W$) and
\begin{equation}
\label{fai} f_{a,i}(u) = 1_{[0, a(i+1)]} (u)\sqrt{a}\int_{(\frac{u}{a}-i)\vee 0} ^{1} \psi(x)K^{H} (a(x+i), u)dx.
\end{equation}
Thus, for all $a>0$ and $i\in \N$,
\begin{multline}
\label{Cpsi}
\hspace{3cm} \mathbb{E} \big (d^{2}(a,i)\big )=\Vert f_{a,i}\Vert ^{2} _{\L^{2}([0,N])}  = a^{2H+1}C_{\psi}(H)\\
\mbox{with }~~~ C_{\psi}(H)=  -\frac 1 2 \, \int_0^1\int_0^1 \psi(x)\psi(x')  | x-x'|^{2H} dx\, dx' ~~~
\end{multline}
(see \cite{Bardet}). This formula is essential for the estimation of $H$ (see Section \ref{Appli}).
Using the product formula (\ref{prod})
$$
I_{1}(f)I_{1}(g)= I_{2} (f\otimes g) + \langle f,g \rangle _{\L^{2}([0,N])}
$$
we obtain
\begin{equation*}
V_{N}(a)= \frac{1}{N_{a}}  \sum_{i=1} ^{N_{a}}\Big ( \frac {I_{2} \big( f_{a,i} ^{\otimes 2} \big)+\Vert f_{a,i}\Vert ^{2} _{\L^{2}([0,N])}} {(\mathbb{E} d(a,i))^2}-1 \Big ) = I_{2} (f_{N}^{(a)})
\end{equation*}
where
\begin{equation}
\label{fn} f_{N}^{(a)}= a^{-2H-1} C_{\psi}(H) ^{-1} \frac{1}{N_{a}}  \sum_{i=1} ^{N_{a}}f_{a,i} ^{\otimes 2} .
\end{equation}
As a consequence the wavelet statistic $V_{N }$  can be expressed as  a multiple Wiener-It\^o integral of order two and therefore the chaos expansion techniques can be applied to its analysis.
\vskip0.5cm

\subsection{A multidimensional Central Limit Theorem satisfied by $(V_{N}(a_i)\big)_{1\leq i \leq \ell}$.}

\vskip0.3cm

When the observed process is a fBm with $H<3/4$, the statistic $V_{N}(a)$ satisfies a central limit theorem. This fact is known and we will not insist on this case. We just recall it to situate it in our context. Since $\mathbb{E} I_{2}^{2}(f)= 2! \Vert f\Vert ^{2} _{\L^{2} [0,N] ^{ \otimes 2 }
}$, if  $(a_i)_{1\leq i \leq \ell}$ is a family of integer numbers such that $a_i=i \, a$ for $i=1,\ldots,\ell$ and $a \in \N^*$, it holds that
\begin{eqnarray*}
\Cov (V_{N}(a_p),V_{N}(a_q))&=& 2! (p \, q\, a^2)^{-2H-1} C_{\psi}(H) ^{-2} \frac{1}{N_{a_p}}
\frac{1}{N_{a_q}}\sum_{j=1}^{N_{a_p}} \sum_{j'=1}^{N_{a_q}}\langle
f_{a_p,j}^{\otimes 2}, f_{a_q,j'}^{\otimes 2} \rangle _{\L^{2}([0,N])^{\otimes 2}}\\
&=& 2(p \, q\, a^2)^{-2H-1} C_{\psi}(H) ^{-2} \frac{1}{N_{a_p}}\frac{1}{N_{a_q}} \sum_{j=1}^{N_{a_p}} \sum_{j'=1}^{N_{a_q}} \langle f_{a_p,j}, f_{a_q,j'} \rangle
_{\L^{2}([0,N])} ^{2}.
\end{eqnarray*}
We know (see e.g.  \cite{Bardet} and \cite{BB}) that:
\begin{eqnarray}
\nonumber \langle f_{a_p,j}, f_{a_q,j'}
\rangle _{\L^{2}([0,N])} &= &\mathbb{E}\left(  d(a_p,j) d(a_q,j')  \right)  \\
\label{covd}& =& -\frac 1 2 \, (p \, q\, a^2)^{1/2} a^{2H}   \int_0^1\int_0^1 \psi(x)\psi(x')  | px-qx'+pj-qj'|^{2H} dx\, dx'.
\end{eqnarray}
and using the  Taylor expansion and the property (\ref{cond-psi}) satisfied by $\psi$,
\begin{multline*}
  \langle f_{a_p,j}, f_{a_q,j'} \rangle _{\L^{2}([0,N])
  }^2  = p\, q \,a^{4H+2}{\cal{O}}\big ( 1+ \vert pj-qj'\vert \big)  ^{4H-4Q}\\
\Longrightarrow ~~
|\Cov (V_{N}(a_p),V_{N}(a_q))|\leq C\,  \frac{1}{N^2_{a_q}} \sum_{j=1}^{N_{a_p}} \sum_{j'=1}^{N_{a_q}} {\cal{O}}\big ( 1+ \vert pj-qj'\vert \big)  ^{4H-4Q}.
\end{multline*}
Consequently, if $Q>1$ and $H\in (0,1)$ or if $Q=1$ and $H \in (0,3/4)$,
\begin{multline}
\label{ell}  N_ {a} \, \Cov (V_{N}(a_p),V_{N}(a_q)) \limiteN \ell_1(p,q,H) ~~~\mbox{with}\\
 \ell_1(p,q,H) = \frac 1   {2 \, d_{pq} \,(p\,q)^{2H-1}}\,  \sum_{k=-\infty}^\infty  \Big (\frac 1 {C_{\psi}(H)} \int_0^1\int_0^1 \psi(x)\psi(x')  | px-qx'+k d_{pq}|^{2H} dx\, dx'\Big )^2,~~~
\end{multline}
where $d_{pq}$ is a constant depending only on $p$ and $q$. Moreover, the following limit theorem holds true.
\begin{theorem}\label{Theo_boncas}
Let $(X_t)_{t\geq 0}$ be a fBm, $V_{N}(a)$ be defined by (\ref{VN}) and $L_1^{(\ell)}(H)=(\ell_1(p,q,H))_{1\leq p,q \leq \ell}$ with $\ell_1(p,q,H)$ defined in (\ref{ell}). Then if $Q>1$ and $H\in (0,1)$ or if $Q=1$ and $H \in (0,3/4)$, for all $a>0$,
\begin{eqnarray}\label{TLCFBM}
\sqrt{ N_{a} } \big (V_{N}(i \, a) \big)_{1\leq i\leq \ell} \limiteloiNa {\cal N}_m\big (0\,, \, L_1^{(\ell)}(H)\big ).
\end{eqnarray}
\end{theorem}
{\bf Proof: } It is well-known in the literature (see e.g. \cite{Bardet}). \qed

\vskip0.3cm

The above central limit theorem does not hold for $Q=1$ and $H\in (\frac{3}{4}, 1)$. We will focus on  this case in the following paragraph.

\vskip0.5cm

\subsection{A noncentral limit theorem satisfied by $V_{N}(a)$.}
We need at this point to define the Rosenblatt process. The Rosenblatt process with time interval $ [0,N]$, denoted in the sequel by $(R^{H}_{t})_{t\in [0,N]}$ appears as a limit in the
so-called {\it noncentral limit theorem } (see \cite{DM} or \cite{Ta1}). It is not a Gaussian process and can be defined
through its representation as a double iterated integral with respect to a standard Wiener process (see \cite{T}). More
precisely, the Rosenblatt process with self-similarity order $H\in (\frac{1}{2}, 1)$
is defined by
\begin{equation}
\label{repZ} R^{H}_{t}=\int _{0}^{t}\int_{0}^{t} L^{H}_{t}(y_1,y_2)dW_ {y_{1}} dW_{y_{2}}
\end{equation}
where $(W_{t}, t \in [0,N])$ is a standard Brownian motion and
\begin{equation}\label{defF}L^{H}_{t}(y_1,y_2)= d_H1_{[0,t]}(y_1)1_{[0,t]}(y_2)\int_{ y_{1} \vee y_{2}
}^{t} \partial_{1}K^{H'} (u,y_{1} ) \partial_{1}K^{H'} (u,y_{2} )du,\end{equation}
with  $K^{H}$ the standard kernel defined in (\ref{defK}) appearing in the Wiener integral representation of the fBm,  its derivatives being defined in (\ref{dK}) and
\begin{eqnarray}\label{dH}
H'=\frac{H+1}{2}\quad \mbox{and}\quad d_H= \frac{1}{H+1} \left( \frac{H}{2(2H-1)} \right) ^{-\frac{1}{2}}.
\end{eqnarray}
The random variable $R^{H}_{1}$ is called a Rosenblatt random variable with self-similarity index $H$. A Rosenblatt process is a process having  stationary increments and
\begin{description}
\item{$\bullet$ }it is $H$-self-similar in the sense that for any
$c>0$, $(R^H _{ct})_{t\geq 0} =  ^{(d)} (c^{H} R^{H}_{t})_{t\geq 0}$, where $ " =^{(d)} "$ means equivalence of all finite dimensional
distributions;
\item{$\bullet$ } $\mathbb{E}(|R^H_t|^p)<\infty$ for any $p>0$, and $R^H$ has the
same variance and covariance as  a standard fractional Brownian motion with parameter $H$.
\item{$\bullet$ } the paths of the Rosenblatt  process are H\"older continuous
of order $\delta <H$.
\end{description}
We obtain the following noncentral limit theorem for the wavelet coefficient of the fBm with $H>\frac{3}{4}$. Define
\begin{eqnarray}\label{l2H}
\ell_{2}(H)&=& \Big (\frac  {2H^2(2H-1)}{4H-3}\Big )^{1/2} \frac{\big (\int_{0}^{1} x\,  \psi(x)\, dx \big )^2}{C_{\psi}(H)}
=\frac{1}{d_{2H-1} }\frac{\big (\int_{0}^{1} x\,  \psi(x)\, dx \big )^2}{C_{\psi}(H)}.
\end{eqnarray}
Then,
\begin{theorem}\label{nclt}
Let $(X_t)_{t\geq 0}$ be a fBm and $V_{N}(a)$ be defined by (\ref{VN}). If $Q=1$ and $\frac{3}{4}<H<1$,
\begin{equation*}
\ell^{-1} _{2}(H) \,  N _{a }^{2-2H}  V_{N}(a) \limiteloiNa R^{H_ {0}}_{1}.
\end{equation*}
where $H_{0}=2H-1$, $\ell_{2}(H)$ is defined by (\ref{l2H}) and $R^{H_{0}}_{1}$ is a Rosenblatt random variable given by $(\ref{repZ})$.
\end{theorem}
{\bf Proof: } With $f_{N}^{(a)}$ defined as in (\ref{fn}), we can write
$$
N_a^{2-2H} V_{N}(a)=N_a^{2-2H} \,  I_{2}\big ( f_{N}^{(a)}\big ).
$$
Using the expression (\ref{fai}) of the kernel $f_{a,i}$, one has
\begin{multline*}
f_{N}^{(a)}(y_{1}, y_{2})= \frac 1 {  a^{2H} C_{\psi }(H)}\,  \frac{1}{N_{a}} \, \sum_{i=1}^{N_{a}}1_{[0, a(i+1)]}(y_{1})  1_{[0, a(i+1)]}(y_{2})\\
\times \int_{(\frac{y_{1}}{a}-i)\vee 0}^{1} \int_{(\frac{y_{2}}{a}-i)\vee 0}^{1}  \psi(x) \psi (z)
K^{H}\left( a(x+i), y_{1}\right)K^{H}\left(  a(z+i), y_{2}\right) dxdz.
\end{multline*}
To show that the sequence $\ell^{-1} _{2}(H)N_{a} ^{2-2H} I_{2} \big(  f_{N} ^{(a)} \big) $ converges in law to the Rosenblatt random variable $R^{2H-1}_1 $ it suffices to show that its cumulants converge to the cumulants of $R^{2H-1}_{1}$ (see e.g. \cite{FT} where it has been proven that the law of a multiple integral of order $2$ is given by its cumulants). On the other hand, we know (see \cite{FT}, \cite{NP}) that the $k$-cumulant of a random variable $I_{2}(f)$ in the second Wiener chaos can be computed as follows
\begin{equation}\label{cum}
c_{k}(I_{2}(f)) = \int_{[0,1] ^{k}} dy_{1} \ldots dy_{k} f(y_{1}, y_{2}) f(y_{2}, y_{3})\ldots f(y_{k-1}, y_{k}) f(y_{k}, y_{1}).
\end{equation}
In particular the $k$th cumulant of the Rosenblatt random variable $R^{2H-1}_{1} $  is given by
\begin{equation*}
c_{k}(R^{2H-1}_{1})=d_{2H-1} ^{k}(H(2H-1))^{k}\int_{[0,1] ^{k}} \big[ |x_{1}-x_{2}|\cdot \ldots \cdot | x_{k-1}-x_{k}|\cdot |x_{k}-x_{1} | \big] ^{2H-2}dx_{1}\ldots dx_{k}.
\end{equation*}
We compute
\begin{multline*}
 c_{k} \big( C_{\psi }(H) N_{a} ^{2-2H} I_{2}  \big(  f_{N} ^{(a)}  \big) \big)\\
 =
  N_{a} ^{(2H-2)k} N_{a}^{-k}\sum_{i_{1}, \ldots , i_{k}=1}^{N_{a}}\int_{[0,1] ^{k}} dy_{1} \ldots dy_{k}\int_{[0,1] ^{2k}}dx_{1}dz_{1} \ldots dx_{k} dz_{k} \psi(x_{1})\psi(z_{1})\psi(x_{2})\psi(z_{2})\ldots \psi(x_{k})\psi(z_{k})\\
   \times K^{H} (a(x_{1} +i_{1} ), y_{1}) K^{H} (a(z_{1} +i_{1} ), y_{2}) K^{H} (a(x_{2} +i_{2} ), y_{2}) K^{H} (a(z_{2} +i_{2} ), y_{3})\times \ldots \\
  \ldots \times K^{H} (a(x_{k-1} +i_{k-1} ), y_{k-1}) K^{H} (a(z_{k} +i_{k} ), y_{k})  K^{H} (a(x_{k} +i_{k} ), y_{k}) K^{H} (a(z_{k} +i_{k} ), y_{1}).
\end{multline*}
Using Fubini Theorem and the fact that (for every $x,x', i, j, a$)
\begin{eqnarray*}
\int_{0}^{a(x+i) \wedge a(x'+j)}K^{H}
(a(x+i), y_{1})K^{H} (a(x'+j), y_{1})dy_{1}&=& R_H(a(x+i),
a(x'+j))\\
&=&a^{2H}R_H(x+i,x'+j)
\end{eqnarray*}
(from the representation of the fBm as a Wiener integral with respect to the Wiener process)
 we get
\begin{eqnarray*}
&&  \hspace{-1cm} c_{k}  \big(C_{\psi }(H) N_{a} ^{2-2H} I_{2}  \big(  f_{N} ^{(a)}  \big) \big)\\
 &=&
  N_{a} ^{(2H-2)k} N_{a}^{-k} \sum_{i_{1}, \ldots , i_{k}=1}^{N_{a}}\int_{[0,1] ^{2k}}dx_{1}dz_{1} \ldots dx_{k} dz_{k} \psi(x_{1})\psi(z_{1})\psi(x_{2})\psi(z_{2})\ldots \psi(x_{k})\psi(z_{k})\\
  && R_H (z_{1}+i_{1}, x_{2}+ i_{2}) R_H (z_{2}+i_{2}, x_{3}+ i_{3})\ldots R_H (z_{k-1}+i_{k-1}, x_{k}+ i_{k})R_H (z_{k}+i_{k}, x_{1}+ i_{1})\\
  &=& N_{a} ^{(2H-2)k}a^{2Hk} N_{a}^{-k} \sum_{i_{1}, \ldots , i_{k}=1}^{N_{a}}\int_{[0,1] ^{2k}}dx_{1}dz_{1} \ldots dx_{k} dz_{k} \psi(x_{1})\psi(z_{1})\psi(x_{2})\psi(z_{2})\ldots \psi(x_{k})\psi(z_{k})\\
 && \times \big [ | z_{1}-x_{2} + i_{1}-i_{2}| \cdot | z_{2}-x_{3} + i_{2}-i_{3} | \ldots \cdot   |z_{k-1} -x_{k}+ i_{k-1}-i_{k} | \cdot | z_{k} -x_{1}+ i_{k}-i_{1}| \big] ^{2H}\\
 &=& N_{a} ^{(2H-2)k} N_{a}^{-k} \sum_{i_{1}, \ldots , i_{k}=1}^{N_{a}}\left( |i_{1}-i_{2} |\ldots \cdot | i_{k-1} -i_{k} | \cdot | i_{k}-i_{1} | \right) ^{2H}\\
 &&\times \int_{[0,1] ^{2k}}dx_{1}dz_{1} \ldots dx_{k} dz_{k} \psi(x_{1})\psi(z_{1})\psi(x_{2})\psi(z_{2})\ldots \psi(x_{k})\psi(z_{k})\left| \left( 1+ \frac{z_{1}-x_{2} }{i_{1}-i_{2}} ) \right) ^{2H} \ldots \left( 1+ \frac{z_{k}-x_{1} }{i_{k}-i_{1}} ) \right) ^{2H}\right| \\
 &\sim & N_{a} ^{(2H-2)k} H^{k} (2H-1) ^{k} N_{a}^{-k} \sum_{i_{1}, \ldots , i_{k}=1}^{N_{a}}\left( |i_{1}-i_{2} |\ldots \cdot | i_{k-1} -i_{k} | \cdot | i_{k}-i_{1} | \right) ^{2H-2}\\
 &&\int_{[0,1] ^{2k}}dx_{1}dz_{1} \ldots dx_{k} dz_{k} \psi(x_{1})\psi(z_{1})\psi(x_{2})\psi(z_{2})\ldots \psi(x_{k})\psi(z_{k})x_{1}z_{1}\ldots x_{k}z_{k}
  \end{eqnarray*}
and we used the fact that the integral of the mother wavelet vanishes and a Taylor expansion of second order of the function $(1+x)^{2H}$. By $a_{n}\sim b_{n}$ we mean that the sequences $a_{n}$ and $b_{n}$ have the same limit as $n\to \infty$. As a consequence, by a  standard Riemann sum argument, the sequence
 \begin{eqnarray*}
N_{a} ^{(2H-2)k}  N_{a}^{-k} \hspace{-4mm}\sum_{i_{1}, \ldots , i_{k}=1}^{N_{a}}\hspace{-3mm} \big( |i_{1}-i_{2} |\times \ldots  \times | i_{k-1} -i_{k} |  | i_{k}-i_{1} | \big) ^{2H-2} = N_{a}^{-k} \hspace{-4mm}\sum_{i_{1}, \ldots , i_{k}=1}^{N_{a}} \hspace{-2mm}\Big( \frac{|i_{1}-i_{2} | \times \ldots \times  | i_{k-1} -i_{k} | | i_{k}-i_{1} | }{N_{a}} \Big) ^{2H-2}
 \end{eqnarray*}
converges to the integral
$$
\int_{[0,1] ^{k}} \big( |x_{1}-x_{2}|\times \ldots \times | x_{k-1}-x_{k}||x_{k}-x_{1} | \big) ^{2H-2}dx_{1}\ldots dx_{k}
$$
and therefore it is clear that the cumulant of $\ell^{-1} _{2}(H)N_{a} ^{2-2H} I_{2}  \big(  f_{N} ^{(a)}  \big) $ converges to the $k$ cumulant of the Rosenblatt random variable $R^{2H-1}_{1}$ (see \cite{Ta1}, \cite{T}). \qed

\vskip0.5cm
In the case of the statistic based on the variations of the fBm, in the case $H\in (3/4, 1)$  the statistic $\displaystyle \frac{1}{N}\sum_{i=0}^{N-1}\frac{(B^{H}_{\frac{i+1}{N}}-B^{H}_{\frac{i}{N}})^{2}}{N^{-2H}}-1$, renormalized by a constant times $N^{2-2H}$, converges in $\L^2(\Omega)$ to a Rosenblatt random variable at time 1 (see \cite{TV}). In the wavelet world, our above result gives only the convergence in law. The following question is then natural: can we get $\L^{2}$-convergence for the renormalized statistic $V_{N}(a)$? The answer is negative and a proof of this fact can be found in the extended version of our paper available on {\tt arXiv}. But we will present below a brief and heuristic argument to see what the $\L^2 $ convergence does not hold. The term $f_{N}^{(a)} $ can be written as
\begin{eqnarray*}
f_{N}^{(a)}(y_{1}, y_{2})&=& \frac 1 {  a^{2H} C_{\psi }(H)}\,  \frac{1}{N_{a}} \, \sum_{i=1}^{N_{a}}  1_{[0, a(i+1)]}(y_{1})  1_{[0, a(i+1)]}(y_{2}) \\
&& \hspace{-2.5cm} \times \Big( 1_{[0,ai] }(y_{1} )1_{[0,ai] }(y_{2} )\int_{0}^{1} \int_{0}^{1} dxdz \psi(x) \psi (z) K^{H}
(a(x+i), y_{1}) K^{H}( a(z+i), y_{2})  \\
&&\hspace{-2cm} +1_{[0,ai] }(y_{1} )1_{[ai, a(1+i)]}(y_{2}) \int_{0}^{1} \int_{\frac{y_{2}}{a}-i}^{1}
 dxdz
\psi(x) \psi (z) K^{H}
(a(x+i), y_{1}) K^{H}( a(z+i), y_{2}) \\
&& \hspace{-2cm} +   1_{[0,ai] }(y_{2} )1_{[ai, a(1+i)]}(y_{1})\int_{\frac{y_{1}}{a}-i}^{1} \int_{0}^{1}
 dxdz
\psi(x) \psi (z) K^{H} (a(x+i), y_{1}) K^{H}( a(z+i), y_{2})  \\
&& \hspace{-1.5cm}+ 1_{[ai, a(1+i)]}(y_{1}) 1_{[ai, a(1+i)]}(y_{2}) \int_{\frac{y_{1}}{a}-i}^{1} \int_{\frac{y_{2}}{a}-i}^{1}
 dxdz
\psi(x) \psi (z) K^{H} (a(x+i), y_{1}) K^{H}( a(z+i), y_{2}) \Big )\\
&& \hspace{-2.5cm} = f_{N}^{(a,1)}(y_{1}, y_{2}) + f_{N}^{(a,2)}(y_{1}, y_{2}) + f_{N}^{(a,3)}(y_{1}, y_{2}) + f_{N}^{(a,4)}(y_{1}, y_{2}).
\end{eqnarray*}
 It can be shown that the terms $N_a^{2-2H} f_{N}^{(a,2)},\, N_a^{2-2H}f_{N}^{(a,3)}$ and $N_a^{2-2H}f_{N}^{(a,4)}$ converge to zero in $\L^{2} ([0,\infty )^{2}) $ as
$N_{a} \to \infty$.  It remains to understand  the convergence of the term $f_{N}^{(a,1)}$. Using again the property (\ref{cond-psi}) of the mother wavelet $\psi$ we can write
\begin{multline*}
f_{N}^{(a,1)}(y_{1}, y_{2})= \frac 1 {  a^{2H} C_{\psi }(H)}\,  \frac{1}{N_{a}}  \sum_{i=0}^{N_{a}}\int_{0}^{1} \int_{0}^{1}dxdz \psi(x) \psi (z)1_{[0, ai]}(y_{1})1_{[0, ai]}(y_{2})\\
\times \left( K^{H}\left( a(x+i), y_{1}\right) -
K^{H}\left( ai, y_{1}\right)\right) \left(
K^{H}\left(
a(z+i), y_{2}\right) - K^{H}\left( ai, y_{2}\right)\right).
\end{multline*}
Therefore, with $\alpha (a,i, x)$ and $\beta (a,i,z)$ located in $[ai,ax +ai]$ and $[ai,az +ai]$ respectively,
\begin{eqnarray*}
I_2 \big ( f_{N}^{(a,1)} \big )
&= &\frac 1 {  a^{2H} C_{\psi }(H)}\,  \frac{1}{N_{a}}  I_2 \Big ( \sum_{i=0}^{N_{a}}\int_{0}^{1} \int_{0}^{1} dxdz \psi(x) \psi (z)1_{[0, ai]}(y_{1})1_{[0, ai]}(y_{2}) \\
&& \hspace{3cm}\times ax \,  \partial_{1}K^{H} (\alpha (a,i,x),
y_{1}) \times az \,  \partial_{1} K^{H}(\beta (a,i,z), y_{2})\Big ).
\end{eqnarray*}
By approximating the points  the points $\alpha (a,i, x)$ and $\beta (a,i,z)$ by $ai$ and since (from an usual  approximation of a sum by a Riemann integral) when $N_a \to \infty$, with $y_1,y_2 \in [0,N]$,
\begin{eqnarray*}
\sum_{i=0}^{N_{a}}1_{[0, ai]}(y_{1})1_{[0, ai]}(y_{2}) \partial_{1}K^{H} (ai, y_{1}) \partial_{1} K^{H}(ai, y_{2})
&\sim & \int_{(y_1 \vee y_2)/a}^{N_a}du \,  \partial_{1}K^{H} (au, y_{1}) \partial_{1} K^{H}(au, y_{2}) \\
&\sim & \frac 1 {d_{2H-1} \,a} \, L^{2H-1}_{N}(y_{1}, y_{2})
\end{eqnarray*}
where $L^{2H-1}_{N}$ is the kernel of the Rosenblatt process with self-similarity index $2H-1$ (see its definition in (\ref{defF})), the sequence $\ell_{2}^{-1}(H)N_{a}^{2-2H} f_{N}^{(a,1)} $ is equivalent (in the sense that it has the same limit pointwise) to $N^{1-2H} L_{N}^{2H-1}$. Then, in some sense, $\ell_{2}^{-1}(H)N_{a}^{2-2H}V_{N}(a)$ is equivalent to $N^{1-2H} I_{2}(L_{N}^{2H-1})=N^{1-2H}R^{2H-1}_{N}= Z_{1}^{2H-1}$ but this equivalence is only in law. The fact  that the sequence $N_{a}^{2-2H}V_{N}(a)$ is not Cauchy in $\L^{2}$ comes from the fact that the sequence $N^{1-2H}R^{2H-1}_{N}$ is not Cauchy in $\L^{2}$ as it can be easily seen by computing the square mean of the difference $N^{1-2H}R^{2H-1}_{N}-M^{1-2H}R^{2H-1}_{M}$.

\section{The Rosenblatt case}\label{SecROS}

\vskip0.2cm

\subsection{Chaotic expansion of the wavelet variation }

We study in this section the limit of the wavelet-based statistic $V_{N}$ given by (\ref{VN}) in the situation when the observed process $X$ is a Rosenblatt process. Throughout this section, assume that $X=R^{H}$ is a Rosenblatt process with self-similarity order $H>1/2$ defined via the stochastic integral representation (\ref{repZ}). We first express the statistic $V_{N}(a)$ in terms of multiple stochastic integrals.  The
wavelet coefficient of the Rosenblatt process can be written as
\begin{eqnarray*}
d(a,i)&=& \sqrt{a} \int_{0}^{1} \psi (x) R^{H} _{a(x+i)} dx \\
&=& \sqrt{a} \int_{0}^{1} \psi(x) dx \Big(    \int_{0}^{a(x+i) }
\int_{0}^{a(x+i)} L^{H}_{a(x+i)}(y_{1}, y_{2}) dW_{y_{1}}dW_{y_{2}}
\Big )\\
&=& I_{2} \left( g_{a,i}(\cdot, \cdot ) \right)
\end{eqnarray*}
where, with  $H'$ and $d_H$ defined in (\ref{dH}),
\begin{multline}
 \nonumber g_{a,i}(y_{1}, y_{2}) =  d_H \,  \sqrt{a} \,  1_{[0,a(i+1)]}(y_{1})1_{[0,a(i+1)]}(y_{2})\\
 \times \int_{\left( \frac{y_{1}\vee y_{2} }{a}-i\right) \vee 0}^1dx \,  \psi (x)\,  \Big(\int _{y_{1} \vee y_{2}}^{a(x+i)}\partial_{1}K^{H'} (u,y_{1} ) \partial_{1}K^{H'} (u,y_{2} )du\Big) \label{gai}
\end{multline}
for every $y_{1}, y_{2} \geq 0$. The product formula for multiple stochastic integrals (\ref{prod}) gives
\begin{equation*}
I_{2}(f)I_{2}(g)= I_{4}(f\otimes g) + 4 I_{2} (f\otimes _{1} g) +
2\langle f,g\rangle _{\L^{2}(0,N] ^{\otimes 2}}
\end{equation*}
if $f,g\in \L^{2} ([0,N]^{2} )$ are two symmetric functions and the
contraction $f\otimes _{1} g$ is defined by
\begin{equation*}
(f\otimes _{1} g)(y_{1},y_{2}) = \int_{0}^{N} f(y_{1},x) g(y_{2},x)
dx.
\end{equation*}
Thus, we obtain
\begin{eqnarray*}
d^{2}(a,i)&=& I_{4} \left( g_{a,i} ^{\otimes 2}\right) + 4 I_{2} \left( g_{a,i} \otimes _{1} g_{a,i} \right)+ 2\Vert g_{a,i}
\Vert ^{2} _{\L^{2}[0,N]^{2}}
\end{eqnarray*}
and note that, since the covariance of the Rosenblatt process is the
same as the covariance of the fractional Brownian motion, we will
also have
$$\mathbb{E} \big (d^2(a,i)\big )= \mathbb{E}\left( I_{2} (g_{a,i} ) \right) ^{2} =2\Vert g_{a,i} \Vert ^{2}
_{\L^{2}[0,N]^{2}}= a^{2H+1}C_{\psi}(H).
$$
Therefore, we obtain the following decomposition for the statistic $V_{N}(a)$:
\begin{eqnarray}
V_{N}(a)&=& a^{-2H-1} C_{\psi }(H) ^{-1} \frac{1}{N_{a}} \left[ \sum _{i=1} ^{N_{a}}I_{4}\left( g_{a,i} ^{\otimes 2} \right)+4
\sum _{i=1} ^{N_{a}}I_{2} \left(
g_{a,i}  \otimes _{1} g_{a,i} \right)\right]=T_{2} + T_{4} \nonumber \\
&&\mbox{with} \qquad \left \{ \begin{array}{ccl} T_2&=&a^{-2H-1} C_{\psi }(H) ^{-1} \frac{4}{N_{a}}\sum _{i=1} ^{N_{a}}  I_{2} \left(
g_{a,i}  \otimes _{1} g_{a,i} \right)\\
T_4&=&a^{-2H-1} C_{\psi }(H) ^{-1} \frac{1}{N_{a}}\sum _{i=1} ^{N_{a}} I_{4}\left( g_{a,i} ^{\otimes 2} \right) \end{array} \right. . \label{decVN}
\end{eqnarray}
To understand the limit of the sequence $V_{N}$ we need to regard the two terms above (note that similar terms
appear in the decomposition of the variation statistic of the Rosenblatt process, see \cite{TV}). In essence, the following
will happen: the term $T_{4}$ which lives in the fourth Wiener chaos keeps some characteristics of the fBm case (since it has
to be renormalized by $\sqrt{N_{a}}$ except in the case $Q=1$ and $H>\frac{3}{4}$ where the normalization is $N_{a} ^{2-2H}$) and its limit will be
Gaussian (except for $Q=1$ and $H>\frac{3}{4}$). Unfortunately, this somehow nice behavior does not affect the limit of $V_{N}$
which is the same as the limit of $T_2$ and therefore it is non-normal (see below).\\
~\\
Now, let us study the asymptotic behavior of  $\mathbb{E} T_{4}^{2}$. From (\ref{decVN}), we have
\begin{equation}
\label{gn}T_{4}= I_{4} (g_{N}^{(a)})\quad \mbox{where}\quad g_{N}^{(a)}=a^{-2H-1} C_{\psi }(H) ^{-1} \frac{1}{N_{a}}
 \sum _{i=1} ^{N_{a}}g_{a,i} ^{\otimes 2},
\end{equation}
and thus, by the isometry of multiple stochastic integrals,
\begin{equation*}
\mathbb{E} T_{4}^{2} = 4! \Vert \tilde{g} _{N}^{(a)} \Vert ^{2} _{\L^{2}[0,N] ^{4} }
\end{equation*}
where by $\tilde{g}_{N}^{(a)} $ we denoted the symmetrization of the function $g_{N}^{(a)} $ is its four variables. Since $\Vert \tilde{g}_{N}^{(a)}\Vert ^{2} _{\L^{2}[0,N] ^{4} }\leq  \Vert g_{N}^{(a)}\Vert ^{2} _{\L^{2}[0,N] ^{4}} $ we obtain
\begin{eqnarray*}
\mathbb{E} T_{4}^{2}&\leq & 4!\,  C_{\psi }(H) ^{-2} a^{-4H-2} \frac{1}{N^2_{a}}  \,  \sum _{i,j=1} ^{N_{a}}\langle g_{a,i}
^{\otimes
2} , g_{a,j} ^{\otimes 2}\rangle _{\L^{2}[0,N] ^{4} } \\
&\leq & 4!\, C_{\psi }(H) ^{-2} a^{-4H-2} \frac{1}{N^2_{a}}  \,\sum _{i,j=1} ^{N_{a}}\langle g_{a,i}, g_{a,j} \rangle
_{\L^{2}(0,N] ^{\otimes 2} } ^{2}.
\end{eqnarray*}
But,
$$\langle g_{a,i}, g_{a,j} \rangle
_{\L^{2}(0,N] ^{\otimes 2} }= \frac{1}{2} \, \mathbb{E} (d(a,i) d(a,j)).
$$
Again, using the fact the fBm and the Rosenblatt process have the same covariance,
we obtain the same behavior (up to a multiplicative constant) as in the case of the fractional Brownian
motion. That is, using (\ref{ell}) and the proof of Theorem \ref{nclt}
\begin{prop} \label{Ros1}
Let $(X_t)_t=(R^{H}_t)_{t\geq 0}$ be a Rosenblatt process and $T_4$ be defined by (\ref{decVN}). If $Q>1$ and $H\in (\frac 1 2 ,1)$ or if $Q=1$ and $H\in (\frac 1 2, \frac{3}{4} )$, then with $\ell_1 (1,1,H)$ defined in (\ref{ell})
\begin{equation} \label{T4cas1}
N_{a}\mathbb{E} (T^2_4) \leq N_{a} 4! \Vert g_{N}^{(a)} \Vert _{\L^{2}[0,N] ^{4} } \limiteN  3 \, \ell_1(1,1,H)
\end{equation}
and, if $Q=1$ and $H\in (\frac{3}{4}, 1)$ then with $\ell_2(H)$ defined in (\ref{l2H}),
\begin{equation}\label{T4cas2}
N_{a}^{4-4H} \, \mathbb{E} (T^2_4) \leq N_{a}^{4-4H} 4! \Vert g_{N}^{(a)} \Vert _{\L^{2}[0,N] ^{4} }\limiteN  3 \, \ell_2(H).
\end{equation}
\end{prop}
The above result gives only an upper bound for the $L^{2}$ norm of the term $T_{4}$; this will be sufficient to obtain the desired limit of the sequence $V_{N}(a)$ in the Rosenblatt case.  It will follow from the following paragraph that  the $L^{2}$ norm of the term $T_{4}$ will be dominated by the $L^{2}$ norm of the term $T_{2}$ and therefore the limit of $T_{2}$ will be the  limit of the statistics $V_{N}(a)$.

\subsection{Asymptotic behavior of the term $T_{2}$}
We study here the term in second Wiener chaos that appears in the decomposition of $V_{N}(a)$. We have
\begin{equation}
\label{t2hn}T_{2}= I_{2}( h_{N}^{(a)})\quad\mbox{with}\quad h_{N}^{(a)}= 4\, \frac 1 {a^{2H+1} C_{\psi}(H)} \frac{1}{N_{a}} \,
 \sum _{i=1} ^{N_{a}}g_{a,i}\otimes _{1} g_{a,i}.
 \end{equation}
We first compute the contraction $g_{a,i}\otimes _{1} g_{a,i}$. We have
\begin{eqnarray*}
&&\hspace{-1cm}(g_{a,i}\otimes _{1} g_{a,i})(y_{1}, y_{2}) = \int_{0}^{N} g_{a,i} (y_{1}, z) g_{a,i} (y_{2}, z) dz\\
& =& a\, d_H ^{2} \,  1_{[0, a(i+1)]}(y_{1} )1_{[0, a(i+1)]}(y_{2} )\int_{0}^{a(i+1)}dz \, \Big [ \int _{\left( \frac{y_{1}\vee z }{a}-i\right) \vee 0}^{1}
dx \, \psi(x) \, \Big( \int_{y_{1}\vee z }^{a(x+i)} \partial_{1}K^{H'} (u,y_{1} )
\partial_{1}K^{H'} (u,z )du\Big)  \Big ] \\
&& \hspace{2cm} \times  \Big [ \int _{\left( \frac{y_{2}\vee z }{a}-i\right) \vee 0}^{1}dx' \,  \psi(x')\, \Big( \int _{ y_{2} \vee z}^{a(x'+i)}\partial_{1} K ^{H'}
(u',y_{2} ) \partial_{1} K ^{H'} (u',z )du'\Big) \Big ]\\
&=&  a\, d_H ^{2} \,  1_{[0, a(i+1)]}(y_{1} )1_{[0, a(i+1)]}(y_{2} ) \Big ( \Big [  \int_{\left(\frac{y_{1}}{a} -i\right) \vee 0} ^{1}dx\, \psi(x) \int_{\left( \frac{y_{2}}{a} -i\right) \vee 0} ^{1}dx'\, \psi(x') \\
&&  \hskip2cm \times  \int_{y_{1}}^{a(x+i)}
 \int _{y_{2} }^{a(x'+i)
}  M(u,y_1,u',y_2) \, du\, du' \int_{0}^{u\wedge u'}  M(u,z,u',z)dz \Big ]
\end{eqnarray*}
where $\displaystyle M(u,y_1,u',y_2)=\partial_{1}K^{H'} (u,y_{1} )\partial_{1} K ^{H'} (u',y_{2} )$
 and $H'=(H+1)/2$. Now, we have already seen that $ \int _{0}^{t\wedge s} K^{H}(t,z) K^{H} (s,z) dz= R_H(t,s)$ with $R_H(t,s)$ given in (\ref{covFBM}) and therefore (see \cite{N}, Chapter 5)
\begin{eqnarray}\label{M}
\int_{0}^{u\wedge u'}  M(u,z,u',z)\,dz =
H'(2H'-1) \, \vert u-u'\vert ^{2H'-2}
\end{eqnarray}
(In fact, this relation can be easily derived from $\int_{0}^{u\wedge
v}K^{H^{\prime}}(u,y_{1})K^{H^{\prime}}(v,y_{1})dy_{1}=R_{H^{\prime}}(u,v)$,
and will be used repeatedly in the sequel). Thus denoting
$\alpha_H=H'(2H'-1)=H(H+1)/2$ and since $\psi$ is $[0,1]$-supported, we obtain
\begin{eqnarray*}
(g_{a,i}\otimes _{1} g_{a,i})(y_{1}, y_{2}) &=& a\, d_H^{2} \alpha_H \,   1_{[0, a(i+1)]}(y_{1} )1_{[0, a(i+1)]}(y_{2} )\int_{\left( \frac{y_{1} }{a}-i\right) \vee 0}^{1}\int_{
\left( \frac{y_{2}}{a}-i\right) \vee 0}^{1}dxdx'\psi(x)\psi(x') \\
&& \hspace{2cm} \times  \int_{y_{1}}^{a(x+i)}
 \int _{y_{2} }^{a(x'+i)
}  \vert u-u'\vert ^{2H'-2} M(u,y_1,u',y_2) \, du\, du'  .
\end{eqnarray*}
By direct computation, it is possible to evaluate the expectation of
$T_{2}^{2}$ as follows (the proof can be found on the extended
version on {\tt arXiv})
\begin{eqnarray}
\nonumber N_{a}^{2-2H}\mathbb{E} T_{2}^{2}\hspace{-3mm}  &\limiteN&
\hspace{-3mm} 32 \,  \frac{\alpha^4_H d_H^{4}}{H(2H-1)C^2_{\psi}(H)
}\Big ( \int_{[0,1] ^{4}}\psi(x) \psi (x') xx'\vert ux-vx'\vert
^{2H'-2} dxdx'dudv\Big )^2 \\
\hspace{-3mm}  &\limiteN&
\hspace{-3mm} C_{T_{2}}^{2}(H)=\frac {32(2H-1)} {H (H+1)^2}\, \Big (\frac { C_{\psi}(H')}{C_{\psi}(H)}\Big)^2= \left( 4d(H) \frac { C_{\psi}(H')}{C_{\psi}(H)}\right) ^{2},\label{ct2}
\end{eqnarray}
where we used  $ \displaystyle \int_{[0,1] ^{2}}\hspace{-4mm}x x'\vert ux-vx'\vert
^{2H'-2} dudv\hspace{-1mm}=\hspace{-1mm}\int_0^x \hspace{-2mm}\int_0^{x'}\hspace{-1mm}\vert u'-v'\vert
^{2H'-2} du'dv'\hspace{-1mm}=\hspace{-1mm}\frac 1 {2H'(2H'-1)} \big (|x|^{2H'}\hspace{-1mm}+\hspace{-0.1mm}|x'|^{2H'}\hspace{-1mm}-\hspace{-0.1mm}|x-x'|^{2H'}\big)$ and $\alpha^4_H\,d^4_H=\frac 1 4 \, (2H-1)^2H^2$.
We do not prove here this estimate because it is a consequence of the
following proposition which show that the sequence
$C_{T_{2}}^{-1}(H)N_{a}^{1-H}T_{2}$ (and therefore the sequence
$V_{N}(a)$) converges in distribution 
to a Rosenblatt random
variable with self-similarity index $H$.

\begin{prop}\label{Ros2}
Let $(R^{H}_{t}) _{t\geq 0 }$ be a Rosenblatt process and let $T_{2}$ be the sequence given by
(\ref{decVN}) and computed from $(R^{H}_{t}) _{t\geq 0 }$. Then, for any  $Q\geq 1$ and $H\in (\frac{1}{2}, 1)$, with $C_{T_{2}}$ given by (\ref{ct2}), there exists a Rosenblatt random variable $R^{H}_{1}$ with self-similarity order $H$ such as
\begin{equation*}
C^{-1}_{T_{2}}(H) \,   N_{a} ^{1-H}  T_{2}\limiteloiN R^{H}_{1}.
\end{equation*}
\end{prop}
{\bf Proof: } This proof follows the lines of the proof of Theorem \ref{nclt}. With $T_{2}= I_{2}\big ( h_{N}^{(a)}\big )$ in mind, as in the proof of Theorem \ref{nclt}, a direct proof that cumulants of the sequence $N_{a}^{1-H}I_{2}\big ( h_{N}^{(a)}\big )$ converge to those of the Rosenblatt process can be given. Indeed, since the random variable $N_{a}^{1-H}I_{2}\big ( h_{N}^{(a)}\big )$ is an element of the second Wiener chaos, its cumulants  can be computed by using formula (\ref{cum}). By using   the key formula (\ref{M})
\begin{eqnarray*}
c_{k} \big(   N_{a}^{1-H}I_{2}\big( h_{N}^{(a)}\big)  \big)&& \\
&& \hspace{-3cm} =N_{a} ^{k(1-H)}\int_{[0,1] ^{k}} dy_{1} \ldots dy_{k}h_{N}^{(a)} (y_{1}, y_{2}) h_{N}^{(a)} (y_{2}, y_{3}) \times \cdots \times h_{N}^{(a)} (y_{k}, y_{1})  \\
&& \hspace{-3cm} =N_{a}^{-Hk} a^{-2Hk} C_{\psi }(H) ^{-k} 4^{k} d_{H}^{2k} \alpha _{H} ^{2k}
\sum_{i_{1}, \ldots , i_{k}=1}^{N_{a}} \int_{[0,1] ^{2k}} \prod_{j=1}^{k} \psi(x_{j} ) \psi(x'_{j} )  dx_{j}dx'_{j} \\
 && \hspace{-3cm}\times  \int_{0}^{a(x_{1}+i_{1})} \hspace{-2mm}\int_{0}^{a(x'_{1}+i_{1})} \hspace{-6mm}  \ldots \hspace{-1mm} \int_{0}^{a(x_{j}+i_{j})}\int_{0}^{a(x'_{j}+i_{j})} \hspace{-2mm} du_{1} du'_{1} \ldots  du_{j} du'_{j}\Big( \prod _{j=1}^{k}\vert u_{j}-u'_{j} \vert \Big) ^{2H'-2} \Big( \prod _{j=1}^{k} \vert u'_{j}-u'_{j+1} \vert \Big) ^{2H'-2}
 \end{eqnarray*}
with the convention $u_{k+1}=u_{1}$. Next, we will make the change of variables $\tilde{u}_{j}= au_{j}$ and $\tilde{u'}_{j}= au'_{j}$ and this will simplify the factors containing $a$. We thus get
\begin{eqnarray*}
c_{k} \big(   N_{a}^{1-H}I_{2}\big( h_{N}^{(a)}\big)  \big)&& \\
&& \hspace{-3cm} = N_{a}^{-Hk} C_{\psi }(H) ^{-k} 4^{k} d_{H}^{2k} \alpha _{H} ^{2k}
\sum_{i_{1}, \ldots , i_{k}=1}^{N_{a}} \int_{[0,1] ^{2k}} \prod_{j=1}^{k} \psi(x_{j} ) \psi(x'_{j} )  dx_{j}dx'_{j} \\
&&\hspace{-3cm} \times  \int_{0}^{x_{1}+i_{1}}\hspace{-2mm}\int_{0}^{x'_{1}+i_{1}} \hspace{-4mm} \ldots \hspace{-2mm} \int_{0}^{x_{k}+i_{k}}\hspace{-2mm} \int_{0}^{x'_{k}+i_{k}}\hspace{-3mm} du_{1} du'_{1} \ldots  du_{k} du'_{k}\Big( \prod _{j=1}^{k}\vert u_{j}-u'_{j} \vert \Big) ^{2H'-2} \Big( \prod _{j=1}^{k} \vert u'_{j}-u'_{j+1} \vert \Big) ^{2H'-2}.
 \end{eqnarray*}
We will note at this point that the vanishing moment property of the wavelet function $\psi $ allows us to replace integration intervals $[0, x_{j}+ i_{j}]$ by intervals $[i_{j},x_{j}+ i_{j}]$. After doing this, we perform the change of variables $\tilde{u}_{j}= u_{j}- i_{j}$, $\tilde{u'}_{j}= u'_{j}- i_{j}$ (for every $j=1,\ldots ,k$) to obtain
 \begin{eqnarray*}
c_{k} \big(   N_{a}^{1-H}I_{2}\big( h_{N}^{(a)}\big)  \big)&& \\
&& \hspace{-3cm} = N_{a}^{-Hk} C_{\psi }(H) ^{-k} 4^{k} d_{H}^{2k} \alpha _{H} ^{2k}
\sum_{i_{1}, \ldots , i_{k}=1}^{N_{a}} \int_{[0,1] ^{2k}} \prod_{j=1}^{k} \psi(x_{j} ) \psi(x'_{j} )  dx_{j}dx'_{j} \\
&&\hspace{-3cm} \times  \int_{0}^{x_{1}}\hspace{-2mm}\int_{0}^{x'_{1}} \hspace{-2mm} \ldots \hspace{-2mm} \int_{0}^{x_{k}}\hspace{-2mm} \int_{0}^{x'_{k}}\hspace{-2mm} du_{1} du'_{1} \ldots  du_{k} du'_{k}\Big( \prod _{j=1}^{k}\vert u_{j}-u'_{j} \vert \Big) ^{2H'-2} \Big( \prod _{j=1}^{k} \vert u'_{j}-u'_{j+1}+i_{j}-i_{j+1}  \vert \Big) ^{2H'-2},
 \end{eqnarray*}
and  then, with changes of variables $\tilde{u}_{j} =\frac{u_{j}}{x_{j}}$ and $\tilde{u'}_{j} =\frac{u'_{j}}{x_{j}}$, we can write
 \begin{eqnarray*}
c_{k} \big(   N_{a}^{1-H}I_{2}\big( h_{N}^{(a)}\big)  \big)&& \\
&& \hspace{-3cm} = N_{a}^{-Hk} C_{\psi }(H) ^{-k} 4^{k} d_{H}^{2k} \alpha _{H} ^{2k}
\sum_{i_{1}, \ldots , i_{k}=1}^{N_{a}} \int_{[0,1] ^{2k}} \prod_{j=1}^{k}x_{j} x'_{j} \, \psi(x_{j} ) \psi(x'_{j} )  dx_{j}dx'_{j} \\
&&\hspace{-3cm} \times  \int_{[0,1] ^{2k}} du_{1}du'_{1} \ldots du_{k}du'_{k} \Big( \prod _{j=1}^{k}\vert x_ju_{j}-x'_ju'_{j} \vert \Big) ^{2H'-2} \Big( \prod _{j=1}^{k} \vert x_ju'_{j}-x'_ju'_{j+1}+i_{j}-i_{j+1}  \vert \Big) ^{2H'-2}\\
&&\hspace{-3cm}= N_{a}^{-Hk} C_{\psi }(H) ^{-k} 4^{k} d_{H}^{2k} \alpha _{H} ^{2k}
\sum_{i_{1}, \ldots , i_{k}=1}^{N_{a}} \int_{[0,1] ^{2k}} \prod_{j=1}^{k}x_{j} x'_{j} \, \psi(x_{j} ) \psi(x'_{j} )  dx_{j}dx'_{j} \\
&&\hspace{-3cm}\times  \sum_{i_{1}, \ldots , i_{k}=1}^{N_{a}} \prod _{j=1} ^{k} \vert i_{j} -i_{j+1} \vert ^{2H'-2}
 \int_{[0,1] ^{2k}} \hspace{-3mm} du_{1}du'_{1}  \ldots du_{k}du'_{k} \Big( \prod_{j=1}^{k}  \vert u_{j}x_{j} -u'_{j}x'_{j} \vert  \Big)^{2H'-2} \prod _{j=1} ^{k}  \Big|    1+ \frac{u'_{j}x'_{j} -u_{j+1}x_{j+1}}{ i_{j} -i_{j+1}}  \Big|  ^{2H'-2}.
\end{eqnarray*}
An analysis of the function $(1+x)^{2H'-2}$ in the vicinity of the origin shows that the cumulant $c_{k} \big(   N_{a}^{1-H}I_{2}\big( h_{N}^{(a)}\big)  \big)$ behaves as (we recall that by $a_{n}\sim b_{n} $ we mean that the sequences $a_{n}$ and $b_{n}$ have the same limit as $n\to \infty$)
\begin{multline*}
c_{k} \big(   N_{a}^{1-H}I_{2}\big( h_{N}^{(a)}\big)  \big) \sim  N_{a}^{-Hk} C_{\psi }(H) ^{-k} 4^{k} d_{H}^{2k} \alpha _{H} ^{2k}
 \int_{[0,1] ^{2k}} \prod_{j=1}^{k} \psi(x_{j} ) \psi(x'_{j} )  dx_{j}dx'_{j} (x_{j}x'_{j})^{k}\\
 \times  \sum_{i_{1}, \ldots , i_{k}=1}^{N_{a}} \prod _{j=1} ^{k} \vert i_{j} -i_{j+1} \vert ^{2H'-2}
 \int_{[0,1] ^{2k}} du_{1}du'_{1} \ldots du_{k}du'_{k} \Big( \prod_{j=1}^{k}  \vert u_{j}x_{j} -u'_{j}x'_{j} \vert \Big)^{2H'-2}.
\end{multline*}
By using an usual Riemann sum convergence, we notice that
\begin{equation*}
N_{a} ^{-Hk} \sum_{i_{1}, \ldots , i_{k}=1}^{N_{a}} \prod _{j=1} ^{k} \vert i_{j} -i_{j+1} \vert ^{2H'-2}=N_{a}^{-k}  \sum_{i_{1}, \ldots , i_{k}=1}^{N_{a}} \prod _{j=1} ^{k} \left( \frac{ \vert i_{j} -i_{j+1} \vert }{ N_{a} } \right) ^{2H'-2}
\end{equation*}
converges as $N_{a} \to \infty$ to the integral
$$
 \int_{[0,1] ^{k} } dx_{1}\ldots dx_{k} \big ( | x_{1} -x_{2}| \cdot | x_{2}-x_{3}| \times \cdots \times  |x_{k}-x_{1}| \big ) ^{2H'-2}
 $$
which is $d_{H} ^{-k} \alpha _{H}^{-k}c_{k} (R^{H} _{1}) $ (here $c_{k}(R^{H}_{1})$ denotes the cumulant of  the random variable $R^{H}_{1}$). We conclude that
 \begin{eqnarray*}
c_{k} \big(   N_{a}^{1-H}I_{2}\big( h_{N}^{(a)}\big)  \big) &\limiteNa &  C_{\psi }(H) ^{-k} 4^{k} d_{H}^{k} c_{k}(R^{H}_{1})  \int_{[0,1] ^{2k}}  \prod_{j=1}^{k}x_{j}x'_{j}  \psi(x_{j} ) \psi(x'_{j} )  dx_{j}dx'_{j} \\
&& \hspace{3cm} \times \int_{[0,1] ^{2k}} du_{1}du'_{1} \ldots du_{k}du'_{k} \Big( \prod_{j=1}^{k}  \vert u_{j}x_{j} -u'_{j}x'_{j} \vert \Big)^{2H'-2}\\
&= &C_{\psi }(H) ^{-k} 4^{k} d_{H}^{k} \alpha _{H} ^{k} \left( \int_{0}^{1} \int_{0}^{1} \psi (x) \psi (x') \frac{1}{2} (x^{2H'} + (x') ^{2H'} -\vert x-x'\vert ^{2H'}dxdx' \right) ^{k} \, c_{k}(R^{H}_{1})\\
&= &C_{\psi }(H) ^{-k} 4^{k} d_{H}^{k}  C_{\psi } (H') ^{k} \, c_{k}(R^{H}_{1})
\end{eqnarray*}
where we recall the notation $C_{\psi }(H) =  - \frac{1}{2} \int_{0}^{1} \int_{0}^{1} \psi (x) \psi (x') \vert x-x'\vert ^{2H'}dxdx'$. As a consequence $C_{T_{2}}^{-1}T_{2} $ converges in law to $R^{H}_{1}$ where $C_{T_{2}}=4d_{H} \frac{C_{\psi }(H') }{C_{\psi }(H)}$. \qed

\vskip0.3cm
\noindent We finally state our main result on the convergence of the wavelet statistic constructed from a Rosenblatt process. Its proof is a consequence of Proposition \ref{Ros1} and \ref{Ros2}.
\begin{theorem}\label{Ros}
Let $(X_t)_{t\geq 0}=(R^{H}_t)_{t\geq 0}$ be a Rosenblatt process and $V_{N}(a)$ be defined by (\ref{VN}). Then, for any  $Q\geq 1$ and $H\in (\frac{1}{2}, 1)$, there exists a Rosenblatt random variable $R^{H}_{1}$ with self-similarity order $H$ such as
\begin{equation*}
C^{-1}_{T_{2}}(H) \,   N_{a}^{1-H}  V_N(a)\limiteloiN R^{H}_{1},
\end{equation*}
where $C_{T_{2}}(H)$ is given by (\ref{ct2}).
\end{theorem}
It is also possible to provide a multidimensional counterpart of
Theorem \ref{Ros} in the case of a vector of scales
$(ai)_{1\leq i\leq \ell}$ where $\ell\in \N^*$.
\begin{theorem}\label{Multros}
Let $(X_t)_{t\geq 0}=(R^{H}_t)_{t\geq 0}$ be a Rosenblatt process and $V_{N}(a)$ be defined by (\ref{VN}). Then for every $Q\geq 1$ and $H\in (\frac{1}{2}, 1)$ it holds
that
\begin{equation*}
\big(  \frac {N_{ai}^{1-H}}{C_{T_{2}}(H)} \,
V_N(ai)\Big) _{1\leq i\leq \ell} \limiteloiN \big(  R^{H}_{1,1},\ldots, R^{H} _{1,\ell}\big)
\end{equation*}
where $ R^{H}_{1,i}$ are normalized Rosenblatt random variables for all
$i=1,\cdots,\ell$ and for all $\lambda _{1},\ldots,\lambda _{\ell} \in \mathbb{R}$ the $k$-th cumulant of the random variable $\sum_{j=1} ^{\ell} \lambda _{j} R^{H}_{1,j}$
is
$$
 \sum_{j_{1},\ldots, j_{k}=1}^{\ell} \lambda _{j_{1}} \ldots \lambda _{j_{k}}c_{k}(R^{H}_{1})$$
where $c_{k}(R^{H}_{1})$ denotes the $k$-th cumulant of the Rosenblatt random variable with self-similarity order $H$.
\end{theorem}
\begin{remark}
Notice that, since the components of the vector $\left(  R^{H}_{1,1},\dots, R^{H} _{1,\ell}\right)$ are random variables in the second Wiener case, its finite dimensional distributions are completely determined by the cumulants. \\
Moreover, we deduce that the asymptotic covariance matrix of the random vector $\big(  N_{a}^{1-H} \,
V_N(ai)\Big) _{1\leq i\leq \ell}$  is
\begin{eqnarray}\label{covRos}
\Sigma_\ell=C^2_{T_{2}}(H) \, \big ( (i\,j)^{1-H} \big) _{1 \leq i,j \leq \ell}.
\end{eqnarray}
\end{remark}
{\bf Proof:} The proof follows the lines of the proof of Theorem \ref{Ros}. Since for every $\lambda _{1},\ldots, \lambda _{m} \in \mathbb{R}$ the linear combination $U_\ell:=\sum_{j=1} ^{\ell} \lambda _{j}\frac {N_{aj}^{1-H}}{C_{T_{2}}(H)} \, V_{N}(aj) $ is a multiple integral of order two, it is possible to compute its cumulants by using the formula (\ref{cum}). We will obtain the following expression for the $k$-th cumulant of $U_\ell$,
\begin{multline*}
\hspace{-4mm} c_{k}(U_\ell)\hspace{-1mm}=\hspace{-1mm}\Big ( \frac {N_{a}^{1-H}}{C_{T_{2}}(H)} \Big )^k \hspace{-2mm}\int_{[0,1] ^{k}}\hspace{-4mm} dy_{1} \ldots dy_{k}
\big( \lambda _{1} h_{N} ^{(a)} (y_{1}, y_{2}) + \ldots + \frac {\lambda _{\ell}}{\ell^{1-H}} h_{N} ^{(\ell a)} (y_{1}, y_{2}) \big) \\ \times \ldots  \times    \big( \lambda _{1} h_{N} ^{(a)} (y_{k}, y_{1}) + \ldots +  \frac {\lambda _{\ell}}{\ell^{1-H}} h_{N} ^{(\ell a)} (y_{k}, y_{1}) \big)
\end{multline*}
and then, as in the proof of Proposition \ref{Ros2}, we can write
\begin{eqnarray*}
c_{k}(U_\ell)&\hspace{-3mm}=&\hspace{-3mm}  \Big ( \frac {4d_{H}^{2} \alpha _{H} ^{2}a^{-2H} N_{a}^{-H}}{ C_{\psi } (H)C_{T_{2}}(H)}\Big )^k \hspace{-4mm} \sum_{j_{1},\ldots , j_{k}=1}^{\ell} \frac {\lambda _{j_{1}}\ldots \lambda _{j_{k}}}{(j_{1}\ldots j_{k})^{H}} \sum_{i_{1}=1} ^{N_{aj_{1}}}\ldots \sum_{i_k=1} ^{N_{aj_{k}}} \int_{[0,1] ^{2k}}\hspace{-2mm} dx_{1}dx'_{1}\ldots dx_{k}dx'_{k}  \prod_{q=1}^k\psi (x_{q}) \psi (x'_{q}) \\
&& \int_{0}^{aj_{1}(x_{1}+i_1)}du_{1} \int_{0}^{aj_{1}(x'_{1}+i_1)}du'_{1} \ldots  \int_{0}^{aj_{k}(x_{k}+i_k)}du_{k} \int_{0}^{aj_{k}(x'_{k}+i_k)}du'_{k}\Big( \prod _{l=1} ^{k} \vert u_{l}-u'_{l}\vert \vert u'_{l}-u_{l+1}\vert \Big)^{2H'-2}
\end{eqnarray*}
where we used again the convention $u_{k+1}=u_{1}$.  We proceed as previously in the proof of Proposition \ref{Ros2} (with changes of variable) and then we obtain
\begin{eqnarray*}
c_{k}(U_\ell)&\hspace{-3mm}=&\hspace{-3mm}  \Big ( \frac {4d_{H}^{2} \alpha _{H} ^{2}\, N_{a}^{-H}}{ C_{\psi } (H)C_{T_{2}}(H)}\Big )^k \hspace{-4mm} \sum_{j_{1},\ldots , j_{k}=1}^{\ell} \frac {\lambda _{j_{1}}\ldots \lambda _{j_{k}}}{(j_{1}\ldots j_{k})^{-1}} \sum_{i_{1}=1} ^{N_{aj_{1}}}\ldots \sum_{i_k=1} ^{N_{aj_{k}}} \int_{[0,1] ^{2k}}\hspace{-2mm} dx_{1}dx'_{1}\ldots dx_{k}dx'_{k}  \prod_{q=1}^k x_{q} x'_{q} \psi (x_{q}) \psi (x'_{q}) \\
&& \int_{[0,1]^{2k} }du_{1} du'_{1} \ldots du_{k} du'_{k}\Big( \prod _{l=1} ^{k} \vert u_{l}x_l-u'_{l}x'_l\vert \vert u'_{l}x'_lj_l-u_{l+1}x_{l+1}j_{l+1}+i_lj_l-i_{l+1}j_{l+1}\vert \Big)^{2H'-2}.
\end{eqnarray*}
Thus, when $N_a \to \infty$,
\begin{eqnarray*}
c_{k}(U_\ell)&\hspace{-3mm} \sim &\hspace{-3mm}  \Big ( \frac {4d_{H}^{2} \alpha _{H} ^{2}\, N_{a}^{-H}}{ C_{\psi } (H)C_{T_{2}}(H)}\Big )^k \hspace{-2mm} \sum_{j_{1},\ldots , j_{k}=1}^{\ell} \frac {\lambda _{j_{1}}\ldots \lambda _{j_{k}}}{(j_{1}\ldots j_{k})^{-1}} \sum_{i_{1}=1} ^{N_{aj_{1}}}\ldots \sum_{i_k=1} ^{N_{aj_{k}}} \int_{[0,1] ^{2k}}dx_{1}dx'_{1}\ldots dx_{k}dx'_{k}  \prod_{q=1}^k x_{q} x'_{q} \psi (x_{q}) \psi (x'_{q}) \\
&&\Big ( \prod _{l=1} ^{k} \vert i_lj_l-i_{l+1}j_{l+1}\vert ^{2H'-2}\Big ) \int_{[0,1]^{2k} }du_{1} du'_{1} \ldots du_{k} du'_{k}\Big( \prod _{l=1} ^{k} \vert u_{l}x_l-u'_{l}x'_l\vert  \Big)^{2H'-2} \\
&\hspace{-3mm} \sim &\hspace{-3mm}  \Big ( \frac {4d_{H}^{2} \alpha _{H} \, a^{2H} \, C_{\psi } (H') \,  N_{a}^{-H}}{ C_{\psi } (H)C_{T_{2}}(H)}\Big )^k \hspace{-2mm} \sum_{j_{1},\ldots , j_{k}=1}^{\ell} \frac {\lambda _{j_{1}}\ldots \lambda _{j_{k}}}{(j_{1}\ldots j_{k})^{-1}} \sum_{i_{1}=1} ^{N_{aj_{1}}}\ldots \sum_{i_k=1} ^{N_{aj_{k}}}  \Big ( \prod _{l=1} ^{k} \vert i_lj_l-i_{l+1}j_{l+1}\vert ^{2H'-2}\Big ) \\
\\
&\hspace{-3mm} \limiteNa  &\hspace{-3mm}  \Big ( \frac {4d_{H}^{2} \alpha _{H} \, a^{2H} \, C_{\psi } (H') \,  }{ C_{\psi } (H)C_{T_{2}}(H)}\Big )^k \hspace{-2mm} \sum_{j_{1},\ldots , j_{k}=1}^{\ell} \lambda _{j_{1}}\ldots \lambda _{j_{k}}\int _{[0,1]^k} \Big ( \prod _{l=1} ^{k} \vert y_l-y_{l+1}\vert ^{2H'-2}\Big )dy_1 \ldots dy_k \\
&\hspace{-3mm} \limiteNa &\hspace{-3mm}  \sum_{j_{1},\ldots , j_{k}=1}^{\ell} \lambda _{j_{1}}\ldots \lambda _{j_{k}} \, c_k(R_1^H),
\end{eqnarray*}
where $R^{H}_{1}$ is a Rosenblatt random variable with index $H$ and $c_{k}(R^{H}_{1}) $ is its $k$-th cumulant. \qed

\vskip0.5cm
\begin{remark}
It is possible and instructive to study the behavior of the term $T_{4}$ in the cases $Q>1$ and $H\in (\frac 1 2 ,1)$ or $Q=1$ and $H\in (\frac 1 2, \frac{3}{4} )$. It can be already seen from its asymptotic variance that it is
very close to the Gaussian case. We can actually show  this term converges in law to a Gaussian random variable. This fact does not influence the limit of the statistic $V_{N}$ but we find that it is interesting from a theoretical point of view.
We will denote by $C_{T_4}(H)$ a
positive constant such that
$$
N_{a} \mathbb{E} T_{4}^{2}\to _{N_{a} \to \infty}C^2_{T_4}(H).
$$
Then the following holds: suppose that $R^{H}$ is a Rosenblatt process with self-similarity order $H$. Suppose that $Q>1$ or $Q=1$ and $H\in (\frac 1 2, \frac{3}{4} )$. Then with $T_4$ defined in (\ref{decVN}),
$$
\sqrt{N_{a}} \ T_{4}  \limiteloiN {\cal N}\big (0\,, \,C^2_{T_{4}}(H) \big ).
$$
The proof can be done by using a criteria in \cite{NOT} in terms of the Malliavin derivatives.
In the other case, {\it i.e.} $Q=1$ and $H\in (\frac{3}{4},1)$, the limit in law of the renormalized $V_{N}$ is a Rosenblatt random
variable.
\end{remark}

\section{Applications and simulations}\label{Appli}
Throughout this section  $X^H=(X^H_t)_{t\in \R _{+}}$ denotes  a fBm with $H\in (0,1)$ or a Rosenblatt process with $H\in (1/2,1)$.

\subsection{Asymptotic normality of the sample variance of approximated wavelet coefficients}
Here a sample $(X^H_0,X^H_1,\cdots,X^H_N)$  of $X^H$ is supposed to be observed. For any couple $(a,b)$, define the following approximations of the wavelet coefficients $d(a,b)$ and of the normalized wavelet coefficients $\tilde d(a,b)$ defined in (\ref{coef}) and (\ref{coefnorm}):
\begin{eqnarray}
e(a,b)=\frac 1 {\sqrt{a}}\, \sum_{k=1}^N X^H_k  \psi \big(\frac k a-b \big )~~\mbox{and}~~ \tilde e(a,b)=\frac {e(a,b)}{a^{H+1/2}C_\psi^{1/2}(H)},
\end{eqnarray}
The above expression  are the usual Riemann approximations. Define also for $a>0$,
\begin{eqnarray}
\widehat V_N(a)= \frac 1 {N_a} \, \sum_{i=1}^{N_a} \big (\tilde e^2(a,i) -1 \big).
\end{eqnarray}
\begin{remark}
These approximations of wavelet coefficients and their sample variance can be directly computed from data for all mother wavelet $\psi$. In the particular case of a multiresolution analysis with orthogonal discrete wavelet transform (that means $a=2^j$), the very fast Mallat's algorithm (similar to FFT for the Fourier's transform) can also be applied to obtain a different approximation $e_M(2^j,k)$ of $d(2^j,k)$. However, when $j\to \infty$, $e(2^j,k)\simeq e_M(2^j,k)$. Indeed, let $\phi$ the scaling function of a multiresolution analysis satisfying $\int \phi(t)dt=1$ and $\int \psi^2(t)dt=1$. Then it can be established that (see for instance \cite{MRT2}) $e_M(2^j,k)=2^{-j/2} \int dt \,  \psi\big(2^{-j}t-k \big)\sum_{s=1}^N \phi(t-s)X_s= 2^{-j/2}\sum_{s=1}^N X_s \int_I dt \,  \psi\big(2^{-j}(t+s)-k \big)   \phi(t)$. When $j\to \infty$, since $\phi$ is compactly supported on $I$, we obtain $ \psi\big(2^{-j}(t+s)-k \big)\sim \psi\big(2^{-j}s-k \big)$ for $s \in I$. Then $e_M(2^j,k) \sim  2^{-j/2}\sum_{s=1}^N X_s \psi\big(2^{-j}s-k \big) \int_I \phi(t)dt=e(2^j,k)$ because $\int \phi(t)dt=1$; more precisely, a first order expansion gives $e_M(2^j,k) = e(2^j,k)+O_P(2^{-j})$, with an approximation error which is negligible with respect to the approximation error computed in the following Lemma \ref{Lemmeapprox}. The very low time consuming of the Mallat's algorithm together with  a straightforward computation of $e(2^j,k)$ without multiresolution analysis, provides a clear advantage of the wavelet based estimator of the parameter $H$ with respect  to the estimators based on a minimization of a criterion (such as maximum likelihood estimators).
\end{remark}
Now the following result can be proved :
\begin{lemma}\label{Lemmeapprox}
Assume that $\psi \in {\cal C}^m(\R)$ with $m \geq 1$ and $\psi$ is $[0,1]$-supported. Let $(a_k)_{k\in \N}$ be a sequence of integer numbers satisfying
$N \, a_N^{-1} \limiteN \infty \quad \mbox{and}\quad
 a_N \limiteN \infty$.
Then, for any $Q\geq 1$,
\begin{eqnarray}\label{diffVN}
\mathbb{E} \big | \widehat V_N(a_N)-V_N(a_N) \big | \leq C \Big ( \frac 1 {\sqrt a_N}+\frac {N^H} { a_N^{H+m}}+\frac {N^{H-Q/2}} {a_N^{(m-Q)/2+H}}1_{(2H-Q> -1)}  \Big ).
\end{eqnarray}
\end{lemma}
{\bf Proof: }  First note that,
\begin{multline*}
\mathbb{E}\big [ (\tilde e(a,i)-\tilde d(a,i))^2 \big ] =\frac {-1} {2\, C_{\psi}(H)} \Big ( \int_0^1 \hspace{-2mm}  \int_0^1 \hspace{-2mm} dtdt' \psi(t)\psi(t')|t-t'|^{2H}+\frac 2 a \int _0^1 \hspace{-2mm} dt \psi(t)  \sum_{k'=0}^{a-1}  \psi \big (\frac {k'}a \big )\Big ( |t+i  |^{2H}-\big |t-\frac {k'} a \big |^{2H}\Big )\\
+\frac 1 {a^2}\sum_{k,k'=0}^{a-1} \psi \big (\frac {k}a \big ) \psi \big (\frac {k'}a \big ) \Big ( \big |i +\frac {k} a  \big |^{2H}+\big |i +\frac {k'} a  \big |^{2H} -\big |\frac {k} a -\frac {k'} a \big |^{2H}\Big ).
\end{multline*}
>From standard Taylor expansion, if $g$ is supposed to be   $m$ times continuously differentiable and  $[0,1]$-supported, 
for all $a>0$,
$$
\Big |\frac 1 a \sum_{k=0}^{a-1}g \big (\frac {k}a \big )-\int_0^1g(t) dt \Big | \leq  \sup_{t\in[0,1]} |g^{(m)}(t)|  \, \frac 1 {a^m}.
$$
and there exists $C$ depending only on $H$, $Q$ and $\psi$ such that $\Big |\int_0^1 \psi(t)(i+t)^{2H}dt \Big | \leq C \big ((1+i)^{2H-Q}\big )$. Therefore, there exists $C$ depending only on $H$, $Q$, $m$ and $\psi$ such that for all $a>0$
$$
\Big |\frac 1 a \sum_{k=0}^{a-1} \psi \big (\frac {k}a \big )\Big | \leq \frac C {a^m}~~\mbox{and}~~ \Big |\frac 1 a \sum_{k=0}^{a-1} \psi \big (\frac {k}a \big )\big | i+\frac k a \big |^{2H}\Big | = C \big ((1+i)^{2H-Q} + \frac {(1+i)^{2H}}{a^m} \big).
$$
Finally, as it was already proved in \cite{Bardet}, there exists $C$ depending only on $H$ and $\psi$ such that for all $m\geq 1$ and $a>0$,
\begin{eqnarray*}
\Big |\frac 1 {a^2}\sum_{k,k'=0}^{a-1} \psi \big (\frac {k}a \big ) \psi \big (\frac {k'}a \big ) \big |\frac {k} a -\frac {k'} a \big |^{2H}-\int_0^1 \hspace{-2mm}  \int_0^1 \hspace{-2mm} dtdt' \psi(t)\psi(t')|t-t'|^{2H} \Big | &\leq &\frac C {a} \\
\Big |\frac 1 a \int _0^1 \hspace{-2mm} dt \psi(t)  \sum_{k'=0}^{a-1}  \psi \big (\frac {k'}a \big )\big |t-\frac {k'} a \big |^{2H} -\int_0^1 \hspace{-2mm}  \int_0^1 \hspace{-2mm} dtdt' \psi(t)\psi(t')|t-t'|^{2H} \Big | &\leq &\frac C {a}.
\end{eqnarray*}
All those inequalities imply that there exists $C$ depending only on $H$, $Q$, $m$ and $\psi$ such that for all $a>0$,
\begin{eqnarray}\label{e-d}
\mathbb{E}\big [ (\tilde e(a,i)-\tilde d(a,i))^2 \big ] \leq C \, \Big ( \frac 1 a + \frac {(1+i)^{2H-Q}}{a^m} + \frac {(1+i)^{2H}}{a^{2m}} \Big ).
\end{eqnarray}
Using  Cauchy-Schwarz's inequality,
\begin{eqnarray*}
\mathbb{E} \big | \widehat V_N(a)-V_N(a) \big | &\leq & \frac 1 {N_a} \, \sum_{i=1}^{N_a} \mathbb{E} \big |\tilde e^2(a,i)-\tilde d^2(a,i) \big |  \\
& \leq & \frac {1} {N_a} \,\sum_{i=1}^{N_a} \sqrt { \mathbb{E} \big [\big (\tilde e(a,i)- \tilde d(a,i) \big )^2 \big ]}\sqrt { \mathbb{E} \big [\big (\tilde e(a,i)+ \tilde d(a,i) \big )^2 \big ]}\\
& \leq & \frac {1} {N_a} \,\sum_{i=1}^{N_a} \sqrt { \mathbb{E} \big [\big (\tilde e(a,i)- \tilde d(a,i) \big )^2 \big ]}\sqrt { \mathbb{E} \big [8\tilde d^2(a,i)+2\big (\tilde e(a,i)- \tilde d(a,i) \big )^2 \big ]} \\
& \leq & \sqrt 2 \, \Big (\frac 1 {N_a} \,\sum_{i=1}^{N_a}  \mathbb{E} \big [\big (\tilde e(a,i)- \tilde d(a,i) \big )^2 \big ] \Big )^{1/2} \Big (\frac 1 {N_a} \,\sum_{i=1}^{N_a} \mathbb{E} \big [4\tilde d^2(a,i)+2\big (\tilde e(a,i)- \tilde d(a,i) \big )^2 \big ] \Big )^{1/2}\\
& \leq & \sqrt 2 \, \Big (\frac 1 {N_a} \,\sum_{i=1}^{N_a}  \mathbb{E} \big [\big (\tilde e(a,i)- \tilde d(a,i) \big )^2 \big ] \Big )^{1/2} \Big ( 4 +2\frac 1 {N_a} \,\sum_{i=1}^{N_a}\mathbb{E} \big [\big (\tilde e(a,i)- \tilde d(a,i) \big )^2 \big ] \Big )^{1/2}
\end{eqnarray*}
since $\mathbb{E} \big [\tilde d^2(a,i)\big ]=1$ by definition. Finally, from inequality (\ref{e-d}) and with $a$ large enough,
\begin{eqnarray*}
\frac 1 {N_a} \,\sum_{i=1}^{N_a}  \mathbb{E} \big [\big (\tilde e(a,i)- \tilde d(a,i) \big )^2 \big ] & \leq &C\, \Big ( \frac 1 a +\frac {\log (N_a)} {N_a  \, a^m}\1_{2H-Q\leq -1}+\frac {N_a^{2H-Q}} {a^m}\1_{2H-Q> -1} + \frac {N_a^{2H}}{a^{2m}}\Big )\\
& \leq &C\, \Big ( \frac 1 a +\frac {\log (N)} {N  \, a^{m-1}}\1_{2H-Q\leq -1}+\frac {N^{2H-Q}} {a^{m+2H-Q}}\1_{2H-Q> -1} + \frac {N^{2H}}{a^{2(m+H)}}\Big ).
\end{eqnarray*}
\qed
~\\
We will use Lemma \ref{Lemmeapprox} to prove the following result.
\begin{prop} \label{Vapprox}
Assume that $\psi \in {\cal C}^m(\R)$ with $m \geq 1$ and $\psi$ is $[0,1]$-supported. Let $(a_k)_{k\in \N}$ be a sequence of integer numbers satisfying
$N \, a_N^{-1} \limiteN \infty \quad \mbox{and}\quad
 a_N \limiteN \infty$. 
Then,
\begin{enumerate}
\item if $X^H$ is a fBm with $0<H<1$ and $Q\geq 2$ or with $0<H<3/4$  and $Q=1$, and if $ N \, a_N^{-1-(1\wedge \frac {2m} 3)}\limiteN 0$, then Theorem \ref{Theo_boncas} holds when $V_N(a)$ is replaced by $\widehat V_N(a)$.
\item if $X^H$ is a fBm with  $3/4<H<1$ and $Q=1$, and if  $ N \, a_N^{-1-(1\wedge \frac {2m} 3)}\limiteN 0$
then  Theorem \ref{nclt} holds when $V_N(a)$ is replaced by $\widehat V_N(a)$.
\item if $X^H$ is a Rosenblatt process with $1/2<H<1$, $Q\geq 1$ and if $ N \, a_N^{-2}\limiteN 0$, then Theorems \ref{Ros} and \ref{Multros} hold when $V_N(a)$ is replaced by $\widehat V_N(a)$.
\end{enumerate}
\end{prop}
{\bf Proof of Proposition \ref{Vapprox}:} The three different cases are respectively obtained from the Markov Inequality. Indeed, for $\varepsilon >0$ and with $\alpha=1/2,~2-2H,~1-H$ (respectively),
$$
\P\big ( N_a^{\alpha}\big | \widehat V_N(a)-V_N(a) \big | >\varepsilon \big )\leq \frac 1 {\varepsilon} \,\frac {N^\alpha}{a^\alpha} \, \mathbb{E} \big | \widehat V_N(a)-V_N(a) \big |.  $$
Using Lemma \ref{Lemmeapprox}, it remains to obtain conditions for insuring
\begin{eqnarray}\label{condTLC}
\frac {N^\alpha}{a^\alpha} \,  \Big ( \frac 1 {\sqrt a_N}+\frac {N^H} { a_N^{H+m}}+\frac {N^{H-Q/2}} {a_N^{\frac {m-Q} 2+H}}\1_{2H-Q> -1}  \Big ) \limiteN 0
\end{eqnarray}
for all $H$ to show the three cases of Proposition \ref{Vapprox}.
\begin{enumerate}
\item For $\alpha=1/2$,  $Q\geq 2$ and $0<H<1$, condition (\ref{condTLC}) with $\alpha=1/2$ leads to
$$\max\big (N\, a_N^{-2}\, , \, N\, a_N^{-1-\frac {2m}{2H+1}}\, , \, N\, a_N^{-1-\frac {m}{2H+1-Q}}\1_{2H-Q>-1} \big )\leq \max\big (N\, a_N^{-2}\, , \, N\, a_N^{-1-\frac {2m}{3}} \big ) \limiteN 0.
$$
It induces the condition $ N \, a_N^{-1-(1\wedge \frac {2m} 3)}\limiteN 0$.\\
For $Q=1$ and any $0<H<3/4$, condition (\ref{condTLC}) leads to
$$\max\big (N\, a_N^{-2}\, , \, N\, a_N^{-1-\frac {2m}{2H+1}}\, , \, N\, a_N^{-1-\frac {m}{2H}} \big ) \limiteN 0.
$$
Since $0<H<3/4$, it also induces the condition $ N \, a_N^{-1-(1\wedge \frac {2m} 3)}\limiteN 0$.
\item if $Q=1$ and $3/4<H<1$, then condition (\ref{condTLC}) with $\alpha=2-2H$ leads to
\begin{multline*}
\max\big (N\, a_N^{-1-\frac 1 {4-4H}}\, , \, N\, a_N^{-1-\frac {2m}{4-2H}}\, , \, N\, a_N^{-1-\frac {m}{3-2H}}\1_{Q=1}\, , \, N\, a_N^{-1-\frac {m}{2-2H}}\1_{Q=2} \big ) \\
\leq \max\big (N\, a_N^{-2}\, , \, N\, a_N^{-1-\frac {2m}{3}} \big ) \limiteN 0,
\end{multline*}
and it also induces the condition $ N \, a_N^{-1-(1\wedge \frac {2m} 3)}\limiteN 0$.
\item if $Q \geq 1$ and $1/2<H<1$, then condition (\ref{condTLC})  with $\alpha=1-H$ leads to
$$
\max\big (N\, a_N^{-1-\frac 1 {2-2H}}\, , \, N\, a_N^{-1-m}\, , \, N\, a_N^{-1-m}\1_{Q=1} \big )\leq \max\big (N\, a_N^{-2}\, , \, N\, a_N^{-1-m} \big ) \limiteN 0,
$$
which leads to the condition $ N \, a_N^{-2}\limiteN 0$ since $m\geq 1$.  \qed
\end{enumerate}

\begin{center}
\begin{figure}
\[
\includegraphics[width=10 cm,height=4.5 cm]{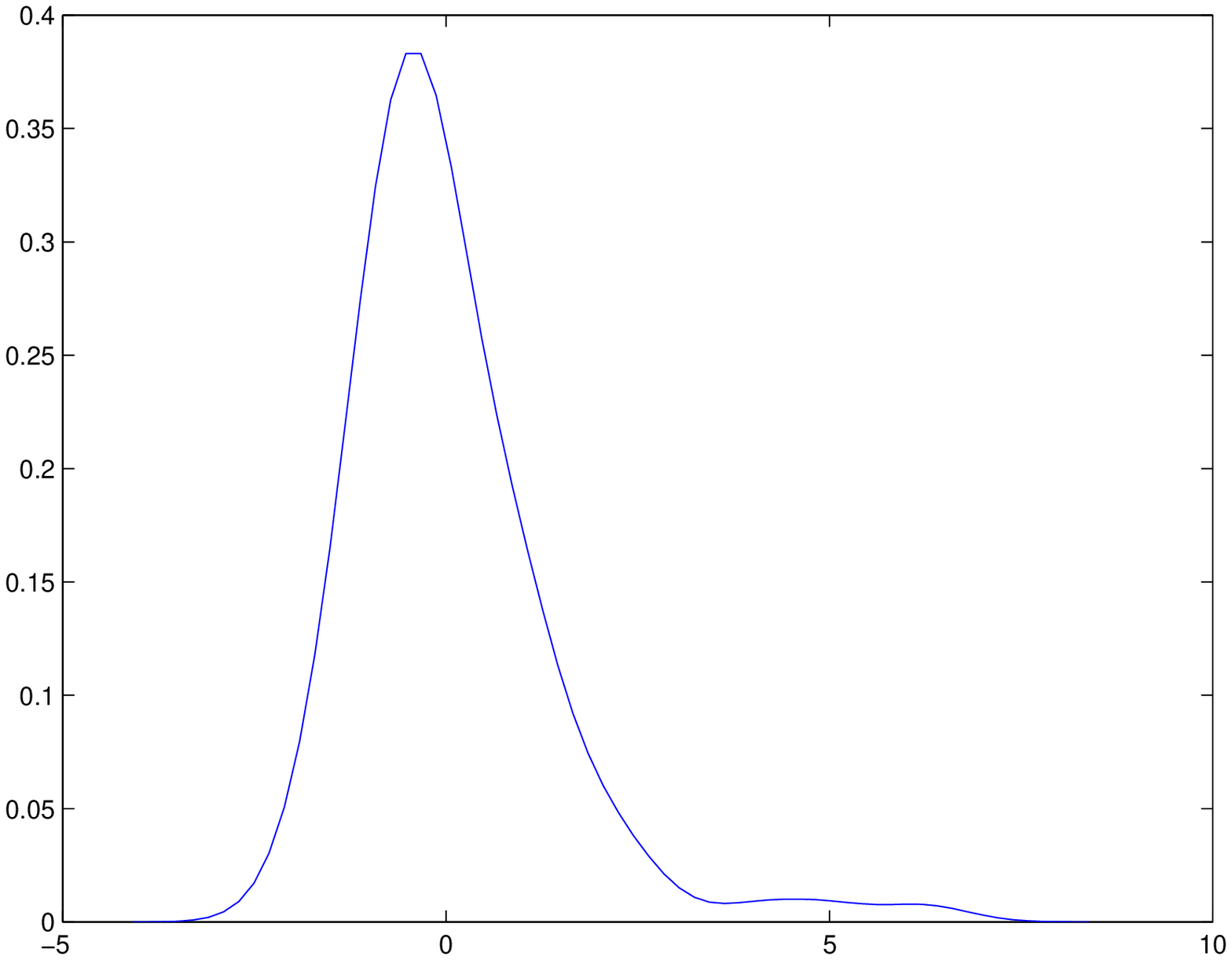}
\]
\caption{\it FFT estimation (Silverman's method) of the density of the limit of $\big (C^{-1}_{T_{2}}(H) \,  \Big (\frac N {a_N} \Big )^{1-H}  \widehat V_N(a_N)\big)_N$ for $H=0.7$, $N=10000$ and $a_N=N^{0.6}$ from $100$ independent replications.}\label{Figure1}
\end{figure}
\end{center}

\begin{figure}
\[
\includegraphics[width=12 cm,height=4.5 cm]{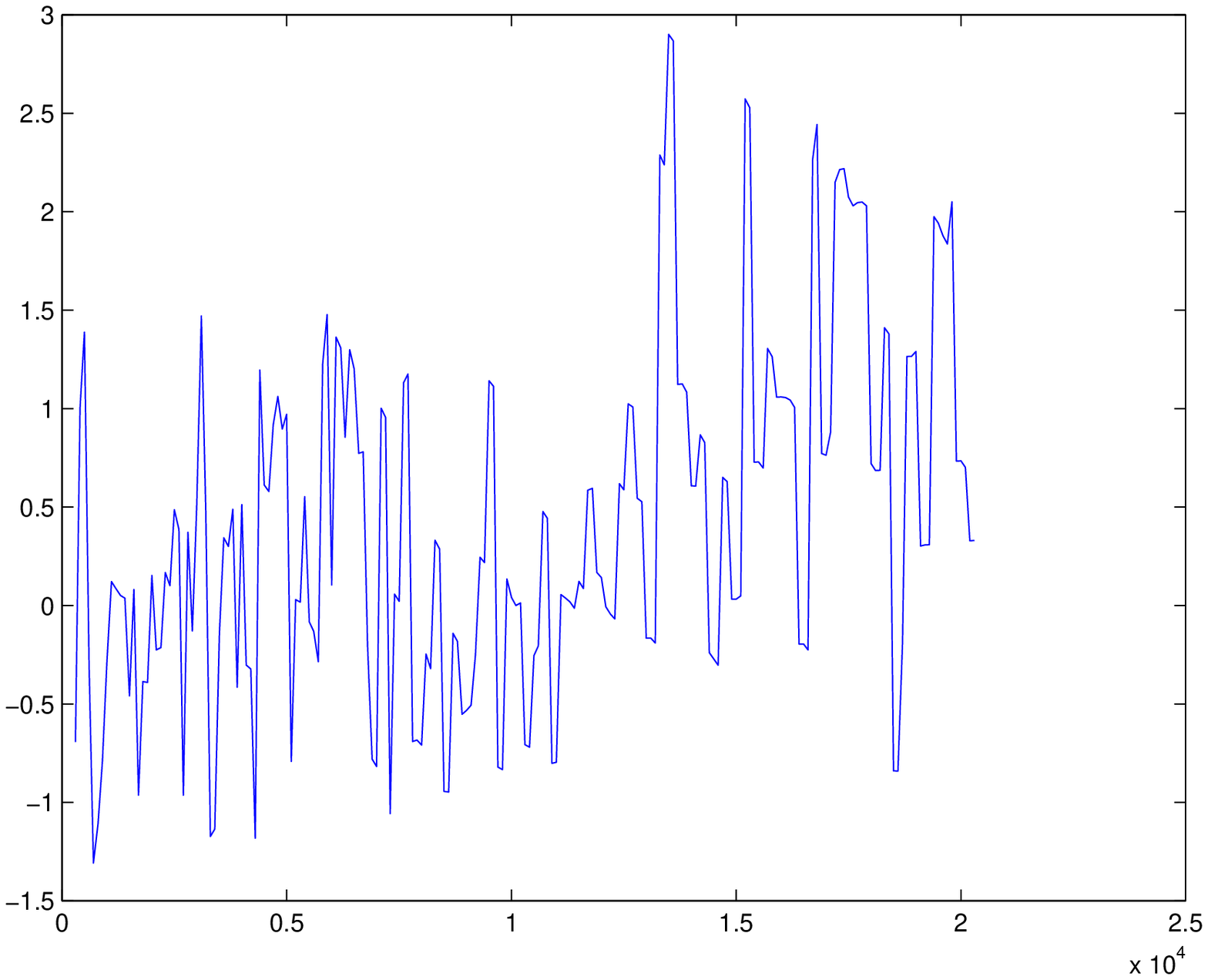}
\]
\caption{\it Convergence of  the sequence $\big (C^{-1}_{T_{2}}(H) \,  \Big (\frac N {a_N} \Big )^{1-H}  \widehat V_N(a_N)\big)_N$ for $H=0.7$ and $a_N=N^{0.5}$.}\label{Figure2}
\end{figure}

\begin{remark}
For a concrete estimation of $H$, the conditions between $N$ and $a_N$ provided in Proposition \ref{Vapprox} do not depend on $H$. Usually, when this conditions depend  on $H$, the convergence rate in the model  can be improved. An adaptive procedure for estimating the smallest order possible for $(a_k)$ could be also built as in the paper \cite{bbj}. Anyway, we do not think that conditions provided on $(a_k)_k$ in Proposition \ref{Vapprox} are optimal. They could be improved by controlling  $\mathbb{E} \big [(\widehat V_N(a)-V_N(a) )^2\big ]$ instead of $\mathbb{E} \big |\widehat V_N(a)-V_N(a) \big |$. However, such computations are very long and tedious in the case of the Rosenblatt process (it requires the computations of fourth-order moments) and we have preferred to avoid them. Moreover, the simulations we give below will show that our results are not so far to be optimal.\\
\end{remark}
Since the case of fBm has been already studied (see for instance \cite{Bardet}) we only provide  below the numerical results when $X^H$ is a Rosenblatt process. Thus, we first exhibit the main result of this paper, {\it i.e.} the limit theorem $C^{-1}_{T_{2}}(H) \,  \Big (\frac N {a_N} \Big )^{1-H}  \widehat V_N(a_N) \limiteloiN R^{H}_{1}$ following the procedure described in the sequel. \\
~\\
{\bf Concrete procedure of simulations:}
\begin{itemize}
\item The samples of Rosenblatt processes are obtained following a similar procedure as the one presented in \cite{pipiras}. With more details, to generate a Rosenblatt process sample $(X^H_j)_{1\leq j \leq N}$:
\begin{enumerate}
\item generate a  sample of length $1+N*m$ (in practice we use $m=100$) of a fBm with parameter $(H+1)/2$. We use a wavelet based method introduced by Sellan (see  \cite{sellan}) with a Daubechies wavelet of order $10$ but a circular matrix embedding method can also be applied (more details  can be found in  \cite{Doukhan}). Note that this sample is normalized and thus $\Var [fBm(1)]=1$. Next, one obtains a sample of length $N*m$ of a fractional Gaussian noise (fGn) defined by the increments of the fBm, {\it i.e.} $fGn(k)=fBm(k+1)-fBm(k)$ with $\Var fGn(k)=1$. 
\item Taqqu proved in \cite{Ta1} that: $\big (\frac 1 {n^H} \sum_{i=1}^{[nt]} (fGn^2(i)-1)\big )_{0\leq t \leq 1} \limiteloin  (R^H_t)_{0\leq t \leq 1}$, where $R^H$ is a (non-normalized) Rosenblatt process, since $(1+H)/2>3/4$ is the parameter of the fGn. Thus, with $n=mN$, $t=j/N$ and $j=1,\ldots,N$, one computes $Y_j= \frac 1 {(mN)^H} \sum_{i=1}^{mj} (fGn^2(i)-1)$. In this way $(Y_j)_{1\leq j\leq N}$ approximatively provides a path of $(R^H_{1/N},R^H_{2/N},\cdots,R^H_1)$.
\item compute now $ X^H_j= \sqrt{\frac{2(2H-1)}{H(H+1)^2}}\,   N^H \, Y_j$ for $j=1,\ldots,N$: using the $H$-self similarity of the Rosenblatt process, $(X^H_1,\ldots,X^H_N)$ is approximatively a trajectory of a normalized ({\it i.e.} $\Var R^H_1=1$) Rosenblatt process.
\end{enumerate}
The Matlab procedures to generate a trajectory of a fBm or a Rosenblatt process can be downloaded from
{\tt http://samos.univ-paris1.fr/-Jean-Marc-Bardet}.
\item Several mother wavelets $\psi$ are used:
\begin{itemize}
\item The Daubechies' wavelet of order $4$, $\psi_4$ (which is such that $Q=4$ and ${\cal C}^1([0,1])$ but not ${\cal C}^2([0,1])$);
\item The Mexican Hat wavelet, $\psi_{MH}$ (which is such that $Q=2$ and $\psi \in {\cal C}^\infty(\R)$ and is essentially compactly supported);
\item The function $\psi_C$ such that $\psi_C(t)=t(t-1)(2t-1)(t^2-t+\frac 1 7)$ for all $t\in [0,1]$  and $\psi_C=0$ elsewhere (which is such that $Q=3$ and $\psi \in {\cal C}^\infty([0,1])$ except in $0$ and $1$).
\end{itemize}
\item The values of constants $C_{\psi}(H)$, $C_{\psi}(H')$ and $C_{T_2}(H)$ are obtained from usual approximations of integrals by Riemann sums.
\end{itemize}
Montecarlo experiments using $100$ independent replications of trajectories are realized for each $H=0.6,\, 0.7,\, 0.8$ and $0.9$ and for $N=500$, $N=2000$ and $N=10000$. The sequence of scales $(a_N)_N$ is selected to be such that $a_N=[N^{0.4}]$, $a_N=[N^{0.5}]$ or $a_N=[N^{0.6}]$.
The following Table \ref{Table1} provides the results of simulations, which are values $\sqrt{\mbox{MSE}}$ of $C^{-1}_{T_{2}}(H) \,  \Big (\frac N {a_N} \Big )^{1-H}  \widehat V_N(a_N)\big)$ for the different choices of $\psi$, $H$, $N$ and $a_N$.
\begin{center}
\begin{table}
\begin{center}
\footnotesize
\begin{tabular}{|c|c|ccc|ccc|ccc|ccc|}
\hline
\multicolumn{2}{|c|}{$H$ } & \multicolumn{3}{|c|}{$0.6$} & \multicolumn{3}{|c|}{$0.7$} & \multicolumn{3}{|c|}{$0.8$} & \multicolumn{3}{|c|}{$0.9$}\\
\cline{1-14}
\multicolumn{2}{|c|}{$a_N$ } & $[N^{0.4}]$& $[N^{0.5}]$& $[N^{0.6}]$ &$[N^{0.4}]$& $[N^{0.5}]$& $[N^{0.6}]$ & $[N^{0.4}]$& $[N^{0.5}]$& $[N^{0.6}]$& $[N^{0.4}]$& $[N^{0.5}]$& $[N^{0.6}]$ \\
\hline
\hline  $N=500$& $\psi_4$ & 3.31&40.0&3.35&1.69& 40.8&1.70&1.08&35.3& 1.49 &1.53&114& 1.05\\
&  $\psi_{MH}$ & 2.23&2.40&2.11&1.51& 1.55&1.80&0.85&0.91 &0.94 &0.70&0.75&0.89 \\
&  $\psi_C$ & 2.29&2.45&2.09&1.76& 1.57&1.49&1.23&0.85&0.93 &1.05&0.74&0.81 \\
\hline  $N=2000$& $\psi_4$ & 6772&3.19&2.99&27935& 7.2&1.64&37951&8.8& 1.56&43888&9.1& 0.75\\
&  $\psi_{MH}$ & 2.09&1.91&1.99&1.33& 1.60&1.40&1.04&1.08&1.30 &0.68&0.67&0.77 \\
&  $\psi_C$ & 1.95&1.81&1.72&1.18& 1.41&1.56&1.21&1.01&1.25 &0.66&0.72& 0.89\\
\hline  $N=10000$&  $\psi_4$ & 2.60&2.77&2.86&1.49& 1.41&1.43&0.67&0.84& 0.93&0.52&0.46& 0.59\\
&  $\psi_{MH}$ & 1.86&1.79&2.07&1.10& 1.24&1.43&0.72&0.70 &0.83 &0.51&0.54&0.60 \\
&  $\psi_C$ & 1.55&1.86&1.72&1.01& 1.08&1.29&0.64&0.73 &0.75 &0.51&0.53&0.56 \\
\hline
\end{tabular}
\end{center}
\caption {\it Computation for different choices of $\psi$, $H$, $N$ and $a_N$ of $\sqrt{\mbox{MSE}}$ of $C^{-1}_{T_{2}}(H) \,  \big (\frac N {a_N} \big )^{1-H}  \widehat V_N(a_N)$ from $100$ independent replications.}\label{Table1}
\end{table}
\end{center}
The main points to highlight from these simulations are the following:
\begin{itemize}
\item globally, for $N$ and $a_N$ large enough then $\big (C^{-1}_{T_{2}}(H) \,  \Big (\frac N {a_N} \Big )^{1-H}  \widehat V_N(a_N)\big)_N$ seems to converge in distribution to a centered distribution with a variance close to $1$;
\item if $a_N$ is not large enough (in these simulations, $a_N=[N^{0.4}]$), there is a bias which clearly appears in the case of $\psi_4$ since this wavelet function is not regular enough. In this way we may see that the conditions required for $(a_N)$ in Proposition \ref{Vapprox} are close to be optimal.
\item Since Rosenblatt processes are generated from an approximation algorithm based on a wavelet synthesis, the larger $H$ the higher the octave required to obtain a convenient trajectory. Because memories of computers are finite and we expect a low consuming time, the generator of Rosenblatt process has a lightly smaller variance than it should have in the case $H=0.9$. Moreover in such a case, there is a  very slow convergence rate, {\it i.e. $(N/a_N)^{0.1}$}, in the limit theorem.
\end{itemize}
An example of the estimation of the limit density is also presented in Figure \ref{Figure1} in the case $H=0.7$, $N=10000$ and $a_N=N^{0.6}$. Such a density is quite similar to a standard Gaussian density but a Kolmogorov-Smirnov test invalides the hypothesis that this distribution is a ${\cal N}(0,1)$ law. This result should be compared with the numerical simulation of the Rosenblatt density given in \cite{ToTu}.\\
Finally, Figure \ref{Figure2} shows the convergence of the sequence $\big (C^{-1}_{T_{2}}(H) \,  \Big (\frac N {a_N} \Big )^{1-H}  \widehat V_N(a_N)\big)_N$ when $N$ increases. This sequence does not seem to converge in $\L^{2}(\Omega)$ as it is claimed in Section \ref{SecFBM}.
\subsection{Estimation of $H$}
Here we consider that a sample $(X_1,\ldots,X_N)$ of $X=\sigma^2\, X^H$ with $X^H$ a fBm or a Rosenblatt process is known and $H$ and $\sigma^2>0$ are unknown. Denote the sample variance of wavelet coefficients
$$
\widehat I_N(a_N):=\frac 1 {N_{a_N}} \sum_{j=1}^{N_{a_N}} e^2(a_N,j).
$$
Then, following Proposition \ref{Vapprox}, one deduces that
\begin{eqnarray}\label{IN}
{N^\alpha_{a_N}}\, \Big (\frac 1 {\sigma^2\, C_{\psi}(H)}\frac {\widehat I_N(ia_N)} {(ia_N)^{2H+1}}-1 \Big )_{1\leq i \leq \ell}\limiteloiN \big (\varepsilon_i)_{1\leq i \leq \ell},
\end{eqnarray}
where the asymptotic distribution $\big (\varepsilon_i)_{1\leq i \leq \ell}$ is a Gaussian distribution (with $\alpha=1/2$ when $X^H$ is a fBm and $Q\geq 2$) or a Rosenblatt distribution as defined in Theorem 4 (with $\alpha=1-H$ when $X^H$ is a Rosenblatt process). Therefore, from the so-called Delta-Method, we also have
\begin{eqnarray}
\label{Delta}
N_{a_N}^{\alpha} \,\Big (\log \big (\widehat I_N(ia_N)\big ) -(2H+1) \log(ia_N)-\log(\sigma^2\, C_{\psi}(H))\Big )_{1\leq i \leq \ell}\limiteloiN \big (\varepsilon_i)_{1\leq i \leq \ell}.
\end{eqnarray}
Therefore, a log-log-regression of $\big ( \widehat I_N(ia_N)\Big )_{1\leq i \leq \ell}$ by $\big (i a_N\big )_{1\leq i \leq \ell}$ provides an estimator of $H$ (this estimation method has been introduced in \cite{flan2}). Such an estimator is defined by
\begin{eqnarray}\label{formulaH}
\widehat H_N:=
\big (\frac 1 2 \, , \, 0\big)'\cdot  \big (Z'_\ell\, Z_\ell\big )^{-1}  Z'_\ell\, \big (\log(\widehat I_N(ia_N))\big )_{1\leq i \leq \ell}-\frac 1 2,
\end{eqnarray}
where $Z_\ell(i,1)=\log i$ and $Z_\ell(i,2)=1$ for all $i=1,\cdots,\ell$.
Then Proposition \ref{Vapprox} implies
\begin{prop}\label{Hconv}
Assume that $\psi \in {\cal C}^m(\R)$ with $m \geq 1$ and $\psi$ is $[0,1]$-supported. Let $(a_k)_{k\in \N}$ be a sequence of integer numbers satisfying
$N \, a_N^{-1} \limiteN \infty \quad \mbox{and}\quad
 a_N \limiteN \infty$. Let $(X_1,\ldots,X_N)$ an observed sample of $X=\sigma^2\, X^H$ where $X^H$ is a fBm or a Rosenblatt process. Then,
\begin{enumerate}
\item if $X^H$ is a fBm with $0<H<1$ and $Q\geq 2$ or with $0<H<3/4$ and $Q=1$, and if $ N \, a_N^{-1-(1\wedge \frac {2m} 3)}\limiteN 0$, then with $\gamma^2(H,\psi,\ell)>0$ defined in (\ref{gamma2}) and depending only on $H$,  $\psi$ and $\ell$,
\begin{eqnarray}\label{TLCHFBM}
\sqrt{\frac N {a_N}} \,  \big (\widehat H_N-H \big) \limiteloiN {\cal N}\big (0 \, , \, \gamma^2(H,\psi,\ell) \big);
\end{eqnarray}

\item if $X^H$ is a Rosenblatt process with $1/2<H<1$, $Q\geq 1$ and if $ N \, a_N^{-2}\limiteN 0$, then
\begin{eqnarray}\label{TLCHRos}
\big (\frac N {a_N}\big)^{1-H} \big (\widehat H_N-H \big ) \limiteloiN L_{H,\psi,\ell}
\end{eqnarray}
where $L_{H,\psi,\ell} $ defined in (\ref{LH}) is a distribution depending only on $H$,  $\psi$ and $\ell$.
\end{enumerate}
\end{prop}
{\bf Proof of Proposition \ref{Hconv}:} From the results of Proposition \ref{Vapprox}, the relation (\ref{IN}) is clear since $\widehat I_N(a)=\widehat V_N(a)+1$ for all $a>0$. Now using the usual multidimensional Delta-method (see for instance \cite{vdv}) with the transformation function $g(x_1,\ldots,x_\ell)=\big (\log(x_1),\ldots,\log(x_\ell)\big )'$ applied to the limit theorem (\ref{IN}), one obtains
$$
N_{a_N}^\alpha \Big (g \Big (\frac 1 {\sigma^2\, C_{\psi}(H)}\Big (\frac {\widehat I_N(ia_N)} {(ia_N)^{2H+1}}\Big)_{1\leq i \leq \ell}\Big ) -g (1,\ldots,1) \Big ) \limiteloiN J_g(1,\ldots,1) \, \big (\varepsilon_i)_{1\leq i \leq \ell},
$$
where $J_g(1,\ldots,1)$ is the Jacobian matrix of $g$ at point $(1,\ldots,1)$. Therefore, since $J_g(1,\ldots,1)$ is the identity matrix, one obtains (\ref{Delta}). Then, $2H+1$ can be estimated from an ordinary least square regression and we obtain:
$$
2\, \widehat H_N+1=\big (1  \, , \, 0\big)'\cdot  \big (Z'_\ell\, Z_\ell\big )^{-1}  Z'_\ell\, \big (\log(\widehat I_N(ia_N))\big )_{1\leq i \leq \ell}=M_\ell\, \big (\log(\widehat I_N(ia_N))\big )_{1\leq i \leq \ell},
$$
with $M_\ell$ a $(1\times \ell)$ matrix. Hence, the formula (\ref{formulaH}) can be deduced. Moreover, (\ref{Delta}) implies
$$
N_{a_N}^{\alpha} \,\Big ((2\, \widehat H_N+1) -(2H+1)\Big )_{1\leq i \leq \ell}\limiteloiN M_\ell \, \big (\varepsilon_i)_{1\leq i \leq \ell}.
$$
Therefore, when $X^H$ is a fBm, since $\big (\varepsilon_i)_{1\leq i \leq \ell}\simloi {\cal N}\big (0\,, \, L_1^{(\ell)}(H) \big ) $ from (\ref{TLCFBM}), we deduce (\ref{TLCHFBM}) with
\begin{eqnarray}\label{gamma2}
\gamma^2(H, \psi, \ell):= \frac 1 4 \, M_\ell \,  L_1^{(\ell)}(H)   \, M'_\ell.
\end{eqnarray}
The same trick implies  (\ref{TLCHRos}) when $X^H$ is a Rosenblatt process and
\begin{eqnarray}\label{LH}
L_{H,\psi,\ell} :=\frac 1 2 \, M_\ell \,\big(  R^{H}_{1,1},\ldots, R^{H} _{1,\ell}\big)'
\end{eqnarray}
with the distribution of $\big(  R^{H}_{1,1},\ldots, R^{H} _{1,\ell}\big)$ defined in Theorem \ref{Multros}.
\qed
\begin{remark}
An estimator of $\sigma^2$ can also be provided by this method, with the same convergence rate.
\end{remark}
\begin{remark}
In Proposition \ref{Hconv} we do not study the case where $X^H$ is a fBm, $Q=1$ and $3/4<H<1$. Indeed, in a statistical framework devoted to the estimation of $H$,
since the choice of $\psi$ is arbitrary, there is no reason to chose $\psi$ with $Q=1$ which gives a  worst convergence rate than $\psi$ with $Q\geq 2$.
\end{remark}
We summarize our results concerning the convergence rate, with $m\geq 2$ and $a_N=N^{1/2+\delta}$ with $\delta>0$ arbitrary small:
\begin{enumerate}
\item if $X^H$ is a fBm and $Q\geq 2$ or $Q=1$ and $0<H<3/4$, the convergence rate of $\widehat H_N$ is $N^{1/4-\delta/2}$;
\item if $X^H$ is a Rosenblatt process and $Q\geq 1$, the convergence rate of $\widehat H_N$ is $N^{(1-H)/2-\delta(1-H)}$.
\end{enumerate}
Such convergence rates are weak in a parametric framework. For instance, applied to the increments of a fBm, the convergence rate of the maximum likelihood or the approximated Whittle maximum likelihood estimator are $N^{1/2}$ (see \cite{dahl89} or \cite{fox-taq1}); as far as we know, there are not such results in the case of the Rosenblatt process. An estimator based on quadratic variations has been  studied in \cite{TV} and a noncentral limit theorem is proved with convergence rate $N^{1-H}$; but such an estimator is almost parametric and can not be applied to processes which are not strictly self-similar. As it was previously recalled in the introduction, the wavelet based estimator is interesting because it can also be applied in a semiparametric framework.
\begin{remark}
In a statistical framework when $X^H$ is a Rosenblatt process, $(X_1,\cdots,X_N)$ is known but $H$ unknown, how to obtain confidence interval from (\ref{TLCHRos})? We propose a three step procedure. The first step is to chose $a_N=[N^{1/2+\delta}]$ with $\delta>0$ arbitrary small (for instance $\delta=0.05$) and compute $\widehat H_N$. The second step consists in the  computation of the quantile $q$ (for instance the $97.5\%$-quantile) of the distribution of $L_{\widehat H_N,\psi,\ell}$ (defined by (\ref{TLCHRos})) from Monte-Carlo simulations of Rosenblatt processes with parameter $H=\widehat H_N$. Indeed, it is possible to replace (\ref{TLCHRos}) by:
$$
\big (\frac N {a_N}\big)^{1-\widehat H_N} \big (\widehat H_N-H \big ) - L_{\widehat H_N,\psi,\ell} \limiteloiN 0,
$$
because $L_{\widehat H_N,\psi,\ell} -L_{H,\psi,\ell} \limiteprobaN 0$ and  $\big (\frac N {a_N}\big)^{1-\widehat H_N} \times \big (\frac N {a_N}\big)^{H-1}=\exp\big [(H-\widehat H_N)\log \big(\frac N {a_N}\big)\big] \limiteprobaN 1$ from (\ref{TLCHRos}) and Slutski's Lemma. Finally, the confidence interval of $H$ will be $\big [-q\, \big (\frac N {a_N}\big)^{\widehat H_N-1} , \,  q\, \big (\frac N {a_N}\big)^{\widehat H_N-1}\big ]$ since the density of $L_{\widehat H_N,\psi,\ell}$ is clearly an even function.
\end{remark}
Since there are not  other asymptotic results than ours on  semiparametric estimators of $H$ for a Rosenblatt process, we can only numerically compare the wavelet estimator with other estimators. We have chosen two usual semiparametric estimators of long memory parameters: the local Whittle estimator (denoted $\widehat H_{Whittle}$) defined in \cite{Rob} and the adaptive global log-periodogram estimator (denoted $\widehat H_{LogPer}$) defined in \cite{MouSou}. Remark that such estimators are applied to the increments of $X$. Note also that both the Matlbab softwares of these estimation procedures can be downloaded from {\tt http://samm.univ-paris1.fr/-Jean-Marc-Bardet}. For completing the numerical study, we have considered two cases of wavelet based estimators:
\begin{itemize}
\item $\widehat H_{Haar}$ where $\psi$ is the Haar mother wavelet defined by $\psi(x)=1$ for $0<x<1/2$ and $\psi(x)=-1$ for $1/2<x<1$ and elsewhere $\psi=0$; moreover, we consider $a_N=[N^{0.5}]$ and $\ell= [N^{0.3}]$.
\item $\widehat H_{\psi_C}$ where $\psi_C$ is wavelet function defined above; moreover, we consider $a_N=[N^{0.4}]$ and $\ell=[N^{0.3}]$.
\end{itemize}
For both these estimators the number $\ell$ of scales is chosen as a sequence depending on $N$ because additional simulations have shown that the variance of these estimators is optimal with such a choice.
The results of simulations are given in Table \ref{Table2}.
\begin{center}
\begin{table}\label{Table2}
\footnotesize
\begin{center}
\begin{tabular}{|c|c|cc|cc|cc|cc|}
\hline
 \multicolumn{2}{|c|}{$H$ } & \multicolumn{2}{|c|}{$0.6$} & \multicolumn{2}{|c|}{$0.7$} & \multicolumn{2}{|c|}{$0.8$} & \multicolumn{2}{|c|}{$0.9$}\\
\cline{3-10}
\multicolumn{2}{|c|}{ } & Mean& $\sqrt{\mbox{MSE}}$ &Mean& $\sqrt{\mbox{MSE}}$& Mean& $\sqrt{\mbox{MSE}}$&Mean & $\sqrt{\mbox{MSE}}$\\
\hline
\hline $N=500$&   $\widehat H_{Haar}$  & 0.38 &0.32&0.43 & 0.35&0.47 &0.37  &0.43 & 0.49  \\
\hline  & $\widehat H_{\psi_C}$&0.56& 0.15& 0.62&0.16&0.72& 0.15& 0.84& 0.10 \\
\hline  &$\widehat H_{Whittle}$ &0.62&0.07& 0.69&0.07&0.77&0.08& 0.87& 0.09 \\
\hline  & $\widehat H_{LogPer}$  &0.58&0.17& 0.67&0.16&0.71&0.18  &0.77 & 0.24 \\
\hline $N=2000$&   $\widehat H_{Haar}$  & 0.50 &0.18&0.60 & 0.18&0.70 &0.16  &0.83 & 0.11 \\
\hline  & $\widehat H_{\psi_C}$&0.59& 0.07& 0.65&0.07&0.75& 0.08&0.85&0.07 \\
\hline  &$\widehat H_{Whittle}$ &0.65&0.06& 0.71&0.04&0.80&0.05& 0.88& 0.05 \\
\hline  & $\widehat H_{LogPer}$  &0.60&0.08& 0.65&0.10&0.75&0.12  & 0.82& 0.13 \\
\hline $N=10000$ &  $\widehat H_{Haar}$  & 0.54 &0.16&0.58 & 0.15&0.66 &0.16  & 0.69& 0.22 \\
\hline  & $\widehat H_{\psi_C}$&0.60& 0.04& 0.66&0.05&0.76& 0.05&0.86&0.04 \\
\hline  &$\widehat H_{Whittle}$ &0.64&0.04& 0.72&0.03&0.79&0.03&0.89&0.03 \\
\hline  & $\widehat H_{LogPer}$  &0.60&0.05& 0.67&0.06&0.74&0.08  &0.85 &0.07 \\
\hline
\end{tabular}
\end{center}
\caption {\it Empirical mean and $\sqrt{\mbox{MSE}}$ of $\widehat H_{Haar}$, $\widehat H_{\phi_C}$, $\widehat H_{Whittle}$, $\widehat H_{LogPer}$ for different choices of $H$ and $N$, from $100$ independent replications.}
\end{table}
\end{center}
The conclusions from Table \ref{Table2} are the following:
\begin{itemize}
\item The local Whittle estimator $\widehat H_{Whittle}$ is clearly the most accurate. The wavelet based estimator $\widehat H_{\psi_C}$ is almost as accurate as $\widehat H_{Whittle}$. The estimator $\widehat H_{LogPer}$ seems to be a little less accurate while $\widehat H_{Haar}$ is not satisfying (it can be explained by the fact that this wavelet function is not continuous in $1/2$).
\item It appears that the convergence rate of the estimators depends on $N$ and all the estimators seem to be consistant.
\item The variance of $\widehat H_{LogPer}$ seems to increase as $N^{1-H}$ when $H$ increases unlike both wavelet estimators while theoretical results say that their variances behave as $N^{1-H}$ (see Proposition \ref{Hconv}). An explanation of this phenomenon is the following: the asymptotic variance $\Var (\widehat H_{\psi})$ of $\widehat H_{\psi_C}$ or $\widehat H_{Haar}$ behaves as the $(2,2)$-component of the matrix $\frac {C^2_{T_2}(H)} 4  \big (Z'_\ell\, Z_\ell\big )^{-1}  Z'_\ell\, \Sigma_\ell \, Z_\ell\,  \big (Z'_\ell\, Z_\ell\big )^{-1}\Big (\frac {a_N} N \Big )^{2-2H}$ where $\Sigma_\ell$ is given in (\ref{covRos}). After computations we obtain that
\begin{eqnarray*}
\Var (\widehat H_{\psi})&=& \frac {C^2_{T_2}(H)} 4  \Big ( \frac {\sum_{i=1}^\ell \log i \times  \sum_{i=1}^\ell i^{1-H}-\ell \sum_{i=1}^\ell i^{1-H} \log i}{\ell \sum_{i=1}^\ell \log^2 i-\big (\sum_{i=1}^\ell \log i\big )^2}\Big )^2  \Big (\frac {a_N} N \Big )^{2-2H} \\
&\sim &\frac {(1-H)^2\,C^2_{T_2}(H)} {4\, H^4} \, \ell^{2H-2} \Big (\frac {a_N} N \Big )^{2-2H}\quad (\ell \to \infty).
\end{eqnarray*}
Therefore the larger $\ell$ the smaller $\Var (\widehat H_{\psi})$. Moreover if $N=10000$, $a_N=[N^{0.4}]$ and $\ell=[N^{0.3}]$, with $\psi=\psi_C$, we obtain for $\Var (\widehat H_{\psi}) \simeq (0.024)^2$,  $(0.028)^2$, $(0.030)^2$ and  $\simeq (0.030)^2$ for respectively $H=0.6,\, 0.7,\, 0.8$ and $0.9$. It explains why the different $\sqrt{\mbox{MSE}}$ do not numerically seem to depend on $H$ even if they theoretically depends on $H$.
\end{itemize}


\section* {References}
\begin{thebibliography}{99}


\bibitem{Abry2} P. Abry, P. Flandrin, M.S. Taqqu and D. Veitch, D. (2003),
\newblock Self-similarity and long-range dependence through the wavelet
lens,
\newblock {In P. Doukhan, G. Oppenheim and M.S. Taqqu editors, {Long-range Dependence: Theory and Applications},
 Birkh{\"a}user.}



\bibitem{pipiras}
P. Abry and V. Pipiras (2006),
\newblock  Wavelet-based synthesis of the Rosenblatt process. Signal Processing, 86, 2326-2339.


\bibitem{Bardet}
{J. M. Bardet (2002), } Statistical study of the wavelet
analysis of fractional Brownian motion.  {IEEE Trans. Inform.
Theory 48,  991-999.}

\bibitem{BB} {J.-M. Bardet and P.R. Bertrand (2008), }A nonparametric estimation of the spectral density of a
continuous-time Gaussian process observed at random times. To appear in Scand. J. Statist.


\bibitem{bbj} J.M. Bardet, H. Bibi, H. and A. Jouini (2008) {Adaptive wavelet based estimator of the memory parameter for stationary Gaussian processes.} \emph{Bernoulli}, 14, 691-724.

\bibitem{BLMS}
{J-M. Bardet, G. Lang, E. Moulines and P. Soulier (2000), }Wavelet estimator of long-range dependent processes. { Stat. Inference Stoch. Process. 4, 85-99.}


\bibitem {Beran}{J. Beran (1994), } Statistics for Long-Memory Processes. {Chapman and Hall. }


\bibitem{BrMa} {P. Breuer and P. Major (1983),} Central limit theorems for
nonlinear functionals of Gaussian fields. J. Multivariate Anal. 13, 425-441.

\bibitem{CNT}
{A.  Chronopoulou, C.A. Tudor and F. Viens (2009),} Application of Malliavin calculus and analysis on Wiener space to long-memory parameter estimation for non-Gaussian processes. C.R.A.S.
Math\'ematiques, 347, 663-666.

\bibitem{CNT1}
{A.  Chronopoulou, C.A. Tudor and F. Viens (2009),} Variations and Hurst index estimation for a Rosenblatt process using longer filters,  Electronic Journal of Statistics, 3, 1393-1435

\bibitem{dahl89} R. Dahlhaus (1989),  Efficient parameter estimation
for self-similar processes.  Ann. Statist. 17,  1749-1766.

\bibitem{DM}
{R.L. Dobrushin and P. Major (1979), }Non-central limit
theorems for non-linear functionals of Gaussian fields. {Z. Wahrsch. Verw. Gebiete, 50,  27-52. }

\bibitem{Doukhan}
P. Doukhan, G. Oppenheim, and M.S. Taqqu (Editors) (2003),
\newblock{Theory and applications of long-range
dependence}, Birkh\"auser.

\bibitem{EM}{P. Embrechts and M. Maejima (2002), } Selfsimilar
processes. {Princeton University Press, Princeton, New York. }

\bibitem{flan2} P. Flandrin (1992), Wavelet analysis and synthesis of
fractional Brownian motion. IEEE Trans. on Inform. Theory,
38, 910-917.

\bibitem{fox-taq1} R. Fox and M.S. Taqqu (1986).  Large-sample properties
of parameter estimates for strongly dependent Gaussian time
series. Ann. Statist. 14, 517-532.

\bibitem{FT} R. Fox and M.S. Taqqu (1987).  Multiple stochastic integrals with dependent integrators.
J. Multivariate Anal. 21, 105-127.

\bibitem{GS} L. Giraitis and D. Surgailis (1990). A central limit theorem for quadratic forms in strongly dependent linear variables and its application to asymptotical normality of Whittle's estimate.  Probab. Theory Related Fields,  86, 87-104.

\bibitem{GT} L. Giraitis and M.S. Taqqu (1999). Whittle estimator for finite-variance non-Gaussian time series with long memory.  Ann. Statist.  27, 178-203.

\bibitem{sellan} Y. Meyer, F. Sellan and M.S. Taqqu (1999). \newblock {Wavelets, generalized white noise and fractional integration: the synthesis of fractional Brownian motion.}  J. Fourier Anal. Appl.  5, 465-494.


\bibitem{MRT2}
E. Moulines, F. Roueff and M.S. Taqqu (2008), Central Limit Theorem for the log-regression wavelet estimation of the memory parameter in the Gaussian semi-parametric context. Fractals, 15, 301–313.

\bibitem{MouSou}
E. Moulines and P. Soulier (2003).
\newblock {Semiparametric spectral estimation for fractional
processes.}
\newblock {In P. Doukhan, G. Oppenheim et M.S. Taqqu editors. Theory and applications of long-range
dependenc, 251--301, Birkh\"auser, Boston.}



\bibitem{N} {D. Nualart (2006),}
Malliavin calculus and related topics, 2nd ed. {Springer.}

\bibitem{NOT} {D. Nualart and S. Ortiz-Latorre (2008),} Central limit
theorems for multiple stochastic integrals and Malliavin calculus. Stochastic Process. Appl.  118, 614-628.


\bibitem{NP} {D. Nualart and G. Peccati (2005)}, Central limit theorems for
sequences of multiple stochastic integrals.  Ann. Probab.  33, 173-193.



\bibitem{Rob}
P.M. Robinson (1995).
\newblock {Gaussian semiparametric estimation of long range dependence},
Ann. Statist., 23, 1630-1661.

\bibitem{RT} F. Roueff and M. S. Taqqu (2009). Asymptotic normality of wavelet estimators of the memory parameter for linear processes. To appear in J. Time Ser. Anal.

\bibitem{ST}{G. Samorodnitsky and M. Taqqu (1994), } Stable
Non-Gaussian random variables. {Chapman and Hall, London. }

\bibitem{Ta1}
{M.  Taqqu (1975),} Weak convergence to the fractional Brownian
motion and to the Rosenblatt process. {Z. Wahrsch. Verw. Gebiete 31,  287-302. }

\bibitem{Taqqu3} {M. Taqqu (1978), } A representation for self-similar processes.  Stochastic Processes and their Applications, 7, 55-64.

\bibitem{ToTu}
{S. Torres and C.A. Tudor (2008), } Donsker type theorem for the Rosenblatt process and a binary market model. Stoch. Anal. Appli. 27, 555-573.

\bibitem{T} {C.A. Tudor (2008), }Analysis of the Rosenblatt process.
ESAIM Probability and Statistics, 12, 230-257.

\bibitem{TV} {C.A. Tudor and F. Viens (2009), } Variations and
estimators for the self-similarity order through Malliavin calculus.  Ann. Probab. 37, 2093-2134.

\bibitem{vdv} {Van der Vaart, A. (1998),} Asymptotic statistics. { Cambridge Series in
Statistical and Probabilistic Mathematics, Cambridge.}

\end{thebibliography}
\end{document}